\documentclass[psamsfonts]{article}

    \usepackage{latexsym}
    \usepackage{amsfonts}
    \usepackage{amsmath}
    \usepackage{amssymb}
    \usepackage{amscd}
    \usepackage{amstext}
    \usepackage{amsthm}

    \title{Classification of homomorphisms and dynamical systems}

 \author{
 Huaxin Lin\\
    Department of Mathematics\\
    University of Oregon\\
    Eugene, Oregon 97403-1222\\
    }
    \date{{\bf Dedicated to George Elliott on his 60th birthday}\\
         }

\begin{document}
    \maketitle

    \newcommand{\CA}{$C^*$-algebra}
    \newcommand{\SCA}{$C^*$-subalgebra}
\newcommand{\aue}{approximate unitary equivalence}
    \newcommand{\ayue}{approximately unitarily equivalent}
    \newcommand{\mops}{mutually orthogonal projections}
    \newcommand{\hm}{homomorphism}
    \newcommand{\pisca}{purely infinite simple \CA}
    \newcommand{\andeqn}{\,\,\,\,\,\, {\rm and} \,\,\,\,\,\,}
    \newcommand{\QED}{\rule{1.5mm}{3mm}}
    \newcommand{\morp}{contractive completely
    positive linear map}
    \newcommand{\asmorp}{asymptotic morphism}
    \newcommand{\arrow}{\rightarrow}
    \newcommand{\tdsum}{\widetilde{\oplus}}
    \newcommand{\pa}{\|}  
    \newcommand{\ep}{\varepsilon}
    \newcommand{\id}{{\rm id}}
    \newcommand{\aueeps}[1]{\stackrel{#1}{\sim}}
    \newcommand{\aeps}[1]{\stackrel{#1}{\approx}}
    \newcommand{\dt}{\delta}
    \newcommand{\yu}{\fang}
    \newcommand{\ca}{{\cal C}_1}
\newcommand{\Ad}{{\rm ad}}
    \newcommand{\tr}{{\rm TR}}
    \newcommand{\N}{{\bf N}}
    \newcommand{\C}{{\bf C}}
    \newcommand{\Aut}{{\rm Aut}}
    \newcommand{\Tand}{\,\,\,\text{and}\,\,\,}
\newcommand{\T}{{\mathbb T}}
\newcommand{\R}{{\mathbb R}}
    \newtheorem{thm}{Theorem}[section]
    \newtheorem{Lem}[thm]{Lemma}
    \newtheorem{Prop}[thm]{Proposition}
    \newtheorem{Def}[thm]{Definition}
    \newtheorem{Cor}[thm]{Corollary}
    \newtheorem{Ex}[thm]{Example}
    \newtheorem{Pro}[thm]{Problem}
    \newtheorem{Remark}[thm]{Remark}
    \newtheorem{NN}[thm]{}
    \renewcommand{\theequation}{e\,\arabic{section}.\arabic{equation}}
    \newcommand{\rforal}{{\rm\,\,\,for\,\,\,all\,\,\,}}
\newcommand{\Z}{{\mathbb Z}}
\newcommand{\Q}{{\mathbb Q}}

    \newcommand{\Ik}{ {\cal I}^{(k)}}
    \newcommand{\Iz}{{\cal I}^{(0)}}
    \newcommand{\Ii}{{\cal I}^{(1)}}
    \newcommand{\Ip}{{\cal I}^{(2)}}

\begin{abstract}
Let $A$ be a unital simple \CA\, with tracial rank zero and $X$ be
a compact metric space. Suppose that $h_1, h_2: C(X)\to A$ are two
unital monomorphisms. We show that $h_1$ and $h_2$ are
approximately unitarily equivalent if and only if
$$
[h_1]=[h_2] \,\,\,{\rm in}\,\,\, KL(C(X),A)\andeqn \tau\circ
h_1(f)=\tau\circ h_2(f)
$$
for every $f\in C(X)$ and every trace $\tau$ of $A.$ Adopting a
theorem of Tomiyama, we introduce a notion of approximate
conjugacy for minimal dynamical systems.
Let $\alpha, \beta: X\to X$ be two minimal
homeomorphisms.  Using the above mentioned result, we show that
two dynamical systems are approximately conjugate in that sense if
and only if a $K$-theoretical condition is satisfied. In the case
that $X$ is the Cantor set, this notion coincides with strong
orbit equivalence of Giordano, Putnam and Skau and the
$K$-theoretical condition is equivalent to saying that the
associate crossed product \CA s are isomorphic.

Another application of the above mentioned result is given for
$C^*$-dynamical systems related to a problem of Kishimoto. Let $A$
be a unital simple AH-algebra with no dimension growth and with
real rank zero, and let $\alpha\in Aut(A).$ We prove that if $\alpha^r$
fixes a large subgroup of $K_0(A)$ and has the tracial Rokhlin
property then $A\rtimes_{\alpha}\Z$ is again a unital simple
AH-algebra with no dimension growth and with real rank zero.

\end{abstract}

\section{Introduction}

Let $X$ be a compact metric space and $\alpha,\beta: X\to X$ be
homeomorphisms. Suppose that $\alpha$ and $\beta$ are minimal,
i.e., neither $\alpha$ nor $\beta$ has non-trivial invariant
closed subsets.  Recall that $\alpha$ and $\beta$ are conjugate if
there exists a homeomorphism $\sigma: X\to X$ such that
$\alpha=\sigma\circ \beta\circ \sigma^{-1}.$ They are flip
conjugate if either $\alpha$ and $\beta$ are conjugate or $\alpha$
and $\beta^{-1}$ are conjugate. It was proved by Jun Tomiyama
(\cite{T})  that $\alpha$ and $\beta$ are flip conjugate if and
only if there exists an isomorphism from
$C(X)\rtimes_{\alpha}{\mathbb Z}$ onto $C(X)\rtimes_{\beta}
{\mathbb Z}$ such that it maps $C(X)$ onto $C(X),$ where
$C(X)\rtimes_{\alpha}{\mathbb Z}$ and $C(X)\rtimes_{\beta}
{\mathbb Z}$ are the associated crossed product \CA s.
It is speculated that \CA\, theory may
help to understand the minimal dynamical systems $(X, \alpha)$ and
$(X, \beta).$ This  is further demonstrated in Giordano, Putnam
and Skau's work (\cite{GPS}) on  minimal Cantor systems. They
show, among other things, that two minimal Cantor systems $(X,
\alpha)$ and $(X, \beta)$ are strong orbit equivalent if and only
if the associated crossed products are isomorphic. It is worth to
noting that under the assumption that $X$ is an infinite set and
$\alpha$ and $\beta$ are minimal, the associated crossed product
\CA s are amenable simple \CA s. Given the recent development in
the classification of amenable simple \CA s, it seems possible to
have some $K$-theoretical description of some useful equivalence
relation for minimal dynamical systems as demonstrated in
Giordano, Putnam and Skau's work. The author has proposed to study
some version of approximate conjugacy for minimal dynamical
systems (see \cite{Lnpro}). The original purpose of this research
is to give  a $K$-theoretical description of approximate conjugacy
in minimal dynamical systems.

Diverging from  dynamical systems, consider \hm s from $C(X),$
the \CA\, of continuous functions on a compact metric space, to a
unital simple \CA\, $A.$ It is fundamentally important in topology
to study \hm s from $C(X)$ to $C(Y).$ In \CA\, theory, it is also
fundamentally important to understand \hm s from $C(X)$ into a
unital \CA. The earliest study of this kind is the classical
Brown-Douglass-Fillmore theory (see \cite{BDF1} and \cite{BDF2}).
The BDF-theory classified  monomorphisms from $C(X)$ into the
Calkin algebra $B(l^2)/{\cal K}(l^2).$ The original motivation was
to classify essentially normal operators. There is no doubt
that BDF-theory plays a crucial role in the development of \CA\,
theory, in particular, in the aspect of \CA\, theory related to
$K$-theory and $KK$-theory. One may notice that the Calkin algebra
is a very special (non-separable) \CA. But it is a simple \CA\,
with real rank zero. It becomes clear that the study of
monomorphisms from $C(X)$ into a separable simple \CA\, with real
rank zero is also very important (see \cite{D}, \cite{Lnrr0}, \cite{Lnjfa},
\cite{LnCX}, \cite{Lnfields},
 \cite{GL1}, \cite{GL2} and \cite{GL3}).
In this paper, we prove
that, under the assumption that
 $A$ is a unital simple \CA\, with tracial rank zero,
then monomorphisms from $C(X)$ into $A$ can be classified up to
approximate unitary equivalence by their $K$-theoretical
information. We believe that this result is potentially very
useful. We will demonstrate this by presenting two applications.

Returning to  minimal dynamical systems, it is known that many
crossed product \CA s associated with minimal dynamical systems
are simple \CA s with tracial rank zero. For example it is known
(see for example Theorem 1.15 of  \cite{GPS}) that
$C(X)\rtimes_{\alpha}{\mathbb Z}$ is a $A\T$-algebra of real rank
zero if $X$ is the Cantor set. It follows (see \cite{EG} and
\cite{Lntaf}) that they have tracial rank zero.  More recently it
is shown by Q. Lin and N.C. Phillips (\cite{LqP}) that crossed
products resulted from minimal diffeomorphisms on a manifold are
inductive limits of subhomogeoneous \CA s. Consequently, by
\cite{Lncre}, they have tracial topological rank zero. Thus by the
classification theorem of \cite{Lnduke} these simple \CA s are
classified by their $K$-theory. The classification of
monomorphisms  from $C(X)$ to those simple \CA s leads to the
notion of approximate (flip) conjugacy in minimal dynamical
systems. Briefly speaking, two minimal dynamical systems $(X,
\alpha)$ and $(X, \beta)$ are approximately $K$-conjugate if there
exist two sequences of homeomorphisms $\sigma_n, \gamma_n: X\to X$
such that
$$
\lim_{n\to\infty}f\circ \sigma_n\circ \beta\circ \sigma_n^{-1}=f\circ \alpha
\andeqn
\lim_{n\to\infty}f\circ \gamma_n\circ \alpha\circ \gamma_n^{-1}=f\circ
\beta
$$
for all $f\in C(X)$ and both $\sigma_n$ and $\gamma_n$ satisfy a
$K$-theoretical constrain. A preliminary result (\cite{LM}) shows
that when $X$ is a Cantor set, $\alpha$ and $\beta$ are
approximately $K$-conjugate if and only if
$C(X)\rtimes_{\alpha}\Z\cong C(X)\rtimes_{\beta}\Z.$ In this
paper, we define a $C^*$-version of approximate flip conjugacy and
use the classification of monomorphisms from $C(X)$ into a unital
simple \CA\, of tracial rank zero to give a $K$-theoretical
condition for two minimal dynamical systems being approximate flip
conjugate in that sense.

We present another related application of the above mentioned
classification of monomorphisms from $C(X).$ We study
$C^*$-dynamical systems $(A, \alpha),$ where $A$ is a unital
simple \CA\, with tracial rank zero and $\alpha\in Aut(A)$ which
satisfies a certain Rokhlin property. The Rokhlin property in
ergodic theory was first adopted into operator algebras in the
context of von Neumann algebras by A. Connes (\cite{Con}). It was
adopted by Herman and Ocneanu (\cite{HO}), then by M. R\o rdam
(\cite{Ro}), A. Kishimoto  and more recently by M. Isumi
(\cite{Iz}) in a much more general context of \CA s. Kishimoto has
studied the problem when a crossed product of a simple
A$\T$-algebra $A$ of real rank zero by an automorphism $\alpha\in
Aut(A)$ is again an A$\T$-algebra of real rank zero. A more
general question is when $A\rtimes_{\alpha}\Z$ is a unital simple
AH-algebra with real rank zero if $A$ is a unital simple
AH-algebra. Given the classification theorem for simple separable
amenable \CA\, with tracial rank zero, a similar question is under
what condition  $A\rtimes_{\alpha}\Z$ has tracial rank zero. In
order to make reasonable sense, one has to assume that $\alpha$ is
sufficiently outer. As proposed by A. Kishimoto (\cite{Ks2}), the
right description of ``sufficiently outer" is that $\alpha$ has a
Rokhlin property. One version of Rokhlin property  was introduced
in \cite{OP} called ``tracial Rokhlin property" which is closely
related to, but slightly weaker than, the so-called approximate
Rokhlin property used in \cite{Ks1}( see \ref{DRok1} below). It
was shown in \cite{OP} that the tracial Rokhlin property occurs
quite often. Tracial cyclic Rokhlin property was introduced in
\cite{LO} which is a stronger Rokhlin property (see \ref{DRokc}
below).
It is shown in \cite{LO} that if $A$ has tracial rank zero and
$\alpha$ has the tracial cyclic Rokhlin property then
$A\rtimes_{\alpha}\Z$ has tracial Rokhlin property. Thus one may
apply classification theorem in \cite{Lnduke} to these simple
crossed products. This leads to the question when an automorphism
$\alpha$ has tracial cyclic Rokhlin property. We will apply the
abovew mentioned results about monomorphisms from $C(X)$
to this problem. Among
other things, we will show that, if $\alpha^r_{*0}|_G={\rm id}_G$
for some subgroup $G\subset K_0(A)$ for which
$\rho_A(G)=\rho_A(K_0(A))$ and $\alpha$ has tracial Rokhlin
property then $\alpha$ has tracial cyclic Rokhlin property. This
result implies that, under the same condition,
$A\rtimes_{\alpha}\Z$ is a unital simple AH-algebra with no
dimension growth and with real rank zero if $A$ is. This solves
the generalized version of the Kishimoto problem.

The paper is organized as follows: In Section 2, we list some
conventions that will be used in this paper. In Section 3, we
present the main results. We first give (in subsection 2.1) the
classification of monomorphisms from $C(X)$ into a unital simple
\CA\, with tracial topological rank zero. We then give two
applications of the theorem. One for minimal dynamical systems
(subsection 2.2) and the other for the $C^*$-dynamical systems and
the Rokhlin property (subsection 2.3). In section 4, we give the
proof of the classification theorem mentioned above and also give
proofs of several approximate version of it. In section 5, we
present the proof of the theorems described in subsection 3.2.
Finally, in section 6, we give the proof of the main results in
subsection 2.3.

{\bf Acknowledgments} This work is partially supported
by a grant from NSF of USA and Zhi-Jiang Professorship
in East China Normal University. The author would like
to acknowledge fruitful conversation with N. C. Phillips
and Hiroyuki Osaka.

\section{Notation}

We will use the following conventions:

(1) Let $A$ be a \CA.
$A_{s.a.}$ is the set of all selfadjoint elements of $A$ and
$A_+$ is the positive cone of $A.$

(2) Let $A$ be a \CA\, and $a\in A.$ We write
${\rm Her}(a)$ for the hereditary \SCA\,  generated by
$a,$ i.e., ${\rm Her}(a)={\overline{aAa}}.$

(3) Let $A$ be a \CA\, and $p, q\in A$ be two projections.
We write $p\sim q,$ and say $p$ is equivalent to $q,$ if there
exists a partial isometry $v\in A$ such that $v^*v=p$ and $vv^*=q.$
If $a\in A_+,$
we write $[p]\le [a]$ if $p\sim q$ for some projection
$q\in {\rm Her}(a).$

(4) Let $A$ be a \CA. We denote by $Aut(A)$ the automorphism group
of $A.$ If $A$ is unital and $u\in A$ is a unitary, we denote by
${\rm ad}\, u$ the inner automorphism defined by ${\rm ad}\,
u(a)=u^*au$ for all $a\in A.$

(5) $T(A)$ is the tracial state space of $A.$ Denote by $\rho_A:
K_0(A)\to Aff(T(A))$ the \hm\, induced by
$\rho_A([p])(\tau)=\tau(p)$ for $\tau\in T(A).$

Furthermore, we also use $\rho_A: A_{s.a}\to Aff(T(A))$ for
the \hm\, defined by $\rho_A(a)(\tau)=\tau(a)$ for $a\in A_{s.a}.$

(6) Let $A$ and $B$ be two \CA s and $\phi, \psi: A\to B$ be two
maps. Let $\ep>0$ and ${\cal F}\subset A$ be a finite subset.
 We write
$$
\phi\approx_{\ep}\psi \,\,\,{\rm on}\,\,\, {\cal F},
$$
if
$$
\|\phi(a)-\psi(a)\|<\ep \rforal a\in {\cal F}.
$$
If $B$ is unital and there is a unitary $u\in B$ such that
$$
\|{\rm ad}\, u\circ \phi(a)-\psi(a)\|<\ep\rforal a\in A,
$$
then we write
$$
\phi \aueeps{u}_{\ep} \psi\,\,\,{\rm on}\,\,\,{\cal F}.
$$

(7) Let  $x\in A,$ $\ep>0$ and ${\cal F}\subset A.$ We write
$x\in_{\ep} {\cal F},$ if ${\rm dist}(x, {\cal F})<\ep,$ or there
is $y\in {\cal F}$ such that $\|x-y\|<\ep.$ Let ${\cal F}$ and
${\cal G}$ be subsets of a \CA\, $A,$ we write ${\cal
F}\subset_{\ep} {\cal G},$ if for every $x\in {\cal F},$
$x\in_{\ep}{\cal G}.$

(8) Let $A$ be a separable amenable \CA. We say that $A$ satisfies
the Universal Coefficient Theorem (UCT), if for any
$\sigma$-unital \CA\, $B,$ one has the following short exact
sequence:
$$
0\to ext_{\Z}(K_{*-1}(A), K_{*}(B))
\to KK^*(A, B)\to Hom(K_*(A), K_*(B))\to 0.
$$
Every \CA\, $A$ in the so-called ``bootstrap" class ${\cal N}$
satisfies the UCT.

Let $G$ and $F$ be abelian groups. Denote by $Pext(G,F)$ the group
of pure group extensions in $ext_{\Z}(G, F),$ i.e., those
extensions of $F$ by $G$ so that every finitely generated subgroup
of $F$ lifts. Denote $KL(A, B)=KK(A,B)/Pext(K_{*-1}(A), K_*(B)).$

(9) Let $C_n$ be a commutative \CA\, with $K_0(C_n)={\mathbb
Z}/n{\mathbb Z}$ and $K_1(C_n)=0.$ Suppose that $A$ is a \CA. Then
$K_i(A, {\mathbb  Z}/k{\mathbb  Z})=K_i(A\otimes C_k).$ Let ${\bf
P}(A)$ be the set of equivalence classes of  projections in
$M_{\infty}(A),$ $M_{\infty}(C(S^1)\otimes A),$
$M_{\infty}((A\otimes C_m{\tilde)})$ and
$M_{\infty}((C(S^1)\otimes A\otimes C_m{\tilde )}).$ We have the
following commutative diagram (\cite{Sc}):
$$
\begin{array}{ccccc}
 K_0(A) &\to &K_0(A, {\mathbb  Z}/k{\mathbb  Z})
&\to& K_1(A)\\
\uparrow_{\bf k} & & & &\downarrow_{\bf k}\\
K_0(A) & \leftarrow & K_1(A, {\mathbb  Z}/k{\mathbb  Z})
& \leftarrow & K_1(A)
 \end{array}
$$
As in \cite{DL}, we use the notation
$$
{\underline K}(A)=
\oplus_{i=0,1, n\in {\mathbb  Z}_+} K_i(A;{\mathbb Z}/n{\mathbb Z}).
$$
 By
$Hom_{\Lambda}({\underline K}(A),{\underline K}(B))$ we mean all
\hm s from ${\underline K}(A)$ to ${\underline K}(B)$ which
respect the direct sum decomposition and the so-called Bockstein
operations (see \cite{DL}). It follows from \cite{DL} that if $A$
satisfies the Universal Coefficient Theorem, then
$Hom_{\Lambda}({\underline K}(A),{\underline K}(B)) \cong
KL(A,B).$

(10) Let $\{A_n\}$ be a sequence of \CA s. Set
$l^{\infty}(\{B_n\})=\prod_{n=1}^{\infty}B_n$ ($C^*$-product of
$\{B_n\}$) and $c_0(\{B_n\})=\oplus_{n=1}^{\infty} B_n$
($C^*$-direct sum). We will use $q_{\infty}(\{B_n\})$ for the
quotient $l^{\infty}(\{B_n\})/c_0(\{B_n\}).$

(11)  Let $A=\lim_{n\to\infty}(A_n,\phi_n),$ where $\phi_n: A_n\to
A_{n+1}$ is the connecting \hm. We denote by $\phi_{n, \infty}:
A_n\to A$ the \hm\, induced by the inductive system. $A$ is said
to be an A$\T$-algebra if each $A_n$ has the form $C(\T)\otimes
F_n,$ for some finite dimensional \SCA\, $F_n.$ $A$ is said to be
an AH-algebra if each $A_n$ has the form $P_nM_{k(n)}(C(X_n))P_n,$
where $X_n$ is a finite CW -complex and $P_n\in M_{k(n)}(C(X_n))$
is a projection. We say that $A$ has no dimension growth if there
is an integer $N$ such that ${\rm dim}X_n\le N.$

(11) Let $A$ and $B$ be two \CA\, and $\phi: A\to B$ be a
\morp. Let $\ep>0$ and ${\cal F}\subset A$ be a subset.
The map  $\phi$ is said to be ${\cal F}$-$\ep$-multiplicative,
if
$$
\|\phi(ab)-\phi(a)\phi(b)\|<\ep\rforal a, b\in {\cal F}.
$$

(12) Let $A$ be a \CA, $\{B_n\}$ be a sequence of \CA s and let
$\phi_n: A\to B_n$ be a sequence of \morp s. We say that
$\{\phi_n\}$ is a sequentially asymptotic morphism if
$$
\lim_{n\to\infty}\|\phi_n(ab)-\phi_n(a)\phi_n(b)\|=0
\rforal a\in A.
$$
Let $\Phi: A\to l^{\infty}(\{ B_n\})$ be defined by
$\Phi(a)=\{\phi_n\}$ and let $\phi=\pi\circ \Phi: A\to
q_{\infty}(\{B_n\}),$ where $\pi: l^{\infty} (\{B_n\})\to
q_{\infty}(\{B_n\})$ is the quotient map. Then $\phi$ is a \hm. In
particular, $[\phi]$ gives an element in $Hom(\underline{K}(A),
\underline{K}(q_{\infty}(\{B_n\})).$ It follows that, for any
finite subset ${\cal P}\subset {\bf P}(A),$ for sufficiently large
$n,$ $[\phi_n]|_{\cal P}$ is well defined partial map to
$\underline{K}(B_n).$

Thus, given any finite subset ${\cal P}\subset {\bf P}(A),$
there is $\dt>0$ and a finite subset ${\cal F}\subset A,$
such that, if $\phi: A\to B$ is a ${\cal F}$-$\dt$-multiplicative
\morp, then
$[\phi]|_{\cal P}$ is well defined (see p. 164 of \cite{Lnctaf} 
and also \cite{DG}).

In what follows, given a finite subset ${\cal P}\subset
{\bf P}(A),$ and $\phi$ is a ${\cal F}$-$\dt$-multiplicative
\morp, when we write $[\phi]|_{{\cal P}}$ we mean it is well defined.

(13)  Let $h: C(X)\to A$ be a \hm\, and $O\subset X$ be
an open subset.
We write
$$
{\rm Her}(h(O))=\{h(f): f(t)=0\, \, \, t\in X\setminus O\}.
$$

(14) Let $A$ be a \CA\, and $\{a_n\}$ be a sequence of elements in
$A.$ We say that $\{a_n\}$ is a {\it central sequence} if
$$
\lim_{n\to\infty}\|a_nx-xa_n\|=0\rforal x\in A.
$$

These conventions will be used throughout the paper without further
explanation.

\section{The main results}

\subsection{Monomorphisms from $C(X)$}

In 1970's, Brown, Douglass and Fillmore (\cite{BDF1}, \cite{BDF2}
and \cite{BDF3}) proved that two unital monomorphisms $h_1$ and
$h_2$ from $C(X)$ into the Calkin algebra, $B(l^2)/{\cal K}(l^2),$
where $B(l^2)$ is the \CA\, of bounded operators and ${\cal
K}(l^2)$ is the \CA\, of compact operators on the Hilbert space
$l^2,$ are unitarily equivalent  if and only if they induce the
same $KK$-element. It should be noted that the Calkin algebra is a
simple \CA\, with real rank zero. M. Dadarlat (\cite{D}) showed
that two monomorphisms from $C(X)$ to a unital purely infinite
simple \CA\, are approximately unitarily equivalent if and only if
they give the same element in $KL(C(X),A).$
 We consider two monomorphisms
$h_1,h_2: C(X)\to A,$ where $A$ is a unital simple \CA\, with
tracial (topological) rank zero. One important previous result was
obtained in \cite{GL2} where we assume that $A$ has a unique
trace.

We recall the definition of tracial (topological) rank of
C*-algebras.

\begin{Def}\label{IDtr}
{\rm Let $A$ be a unital simple \CA~ and $k \in {\mathbb N}$. Then
$A$ is said to  have
{\it tracial (topological) rank zero} if and only if
for any finite set ${\mathcal F} \subset A$, and $\ep > 0$ and
any non-zero positive element $a \in A$, there exists a
finite dimensional \SCA\, $B \subset A$
with
${\rm id}_B = p$ such that
\begin{itemize}
\item[$(1)$]
$\|[x, p]\| < \ep$ for all $x \in {\mathcal F}$,
\item[$(2)$] $pxp \in_\ep B$ for all $x \in {\mathcal F}$,
\item[$(3)$]
$[1 - p] \leq [a]$.
\end{itemize}
We write $\tr(A)=0$ if $A$ has tracial (topological) rank zero. }
\end{Def}

Recall that a \CA\, $A$ is said to have the Fundamental Comparison
Property, if, for any two projections $p,\, q\in A,$
$\tau(p)<\tau(q)$ for all $\tau\in T(A)$ implies that $p\sim q'\le
q.$ If $A$ has the Fundamental Comparison Property, then condition
(3) above can be replaced by (3')\,\,\, $\tau(1-p)<\ep$ for all
$\tau\in T(A).$

It is proved in \cite{Lntr0} that if $\tr(A)=0,$ then
$A$ has real rank zero, stable rank one and weakly unperforated
$K_0(A).$ Every simple AH-algebra with slow dimension growth and
with real rank zero has tracial rank zero.
Other simple \CA s that are inductive limits of type I are also
proved to have tracial rank zero (see \cite{Lncre}).
Separable simple
amenable \CA s with tracial rank zero which satisfy the UCT
are classified by their $K$-theory (see \cite{Lnduke} and \cite{Lnann}).

\begin{Def}\label{IDappu}
Let $A$ be a \CA\, and $B$ be a unital \CA. Suppose that $h_1,
h_2: A\to B$ are two maps. We say $h_1$ and $h_2$ are
approximately unitarily equivalent if
$$
h_1\aueeps{u}_{\ep} h_2 \,\,\,{\rm on}\,\,\,{\cal F}
$$
for every finite subset ${\cal F}$ and $\ep>0.$
In other words, there exists a sequence of unitaries
$u_n\in B$ such that
$$
\lim_{n\to\infty}\|{\rm ad}\, u_n\circ h_1(a)-h_2(a)\|=0\rforal a\in A.
$$
\end{Def}

Suppose that $A$ has tracial states and both $h_1$ and $h_2$ are
\hm s. Let $\tau\in T(A).$ It is clear that if $h_1$ and $h_2$ are
approximately unitarily equivalent then $\tau\circ h_1=\tau\circ
h_2.$ In other words, $\tau\circ h_1$ and $\tau\circ h_2$ induce
the same Borel measure on $X.$

The first main result of this paper is the following:

\begin{thm}\label{MT02}
Let $X$ be a compact metric space and $A$ be a unital separable
simple \CA\, with $TR(A)=0.$ Suppose that $h_1: C(X)\to A$ is a
unital monomorphism. For any $\ep>0$ and any finite
subset ${\cal F}\subset C(X),$ there exist $\dt>0$ and a finite
subset ${\cal G}\subset C(X)$ satisfying the following: if
$h_2: C(X)\to A$ is another unital monomorphism such that
$$
[h_1]=[h_2]\,\,\,{\rm in}\,\,\, KL(C(X), A)\andeqn
 |\tau\circ
h_1(f)-\tau\circ h_2(f)|<\dt
$$
for all $f\in {\cal G}$ and $\tau\in T(A),$  then there exists a unitary $u\in A$ such
that
$$
h_1\aueeps{u}_{\ep} h_2\,\,\,{\rm on}\,\,\, {\cal F}.
$$
\end{thm}

Consequently, we have

\begin{thm}\label{MT1}
Let $X$ be a compact metric space and $A$ be a unital
separable simple \CA\, with $TR(A)=0.$
Suppose that $h_1, h_2: C(X)\to A$ are two unital monomorphisms.
Then $h_1$ and $h_2$ are approximately unitarily
equivalent if and only if
$$
[h_1]=[h_2]\,\,\,{\rm in}\,\,\, KL(C(X), A)
$$
and $\tau\circ h_1(a)=\tau\circ h_2(a)$ for all
$a\in C(X).$
\end{thm}

A special case when $A$ has a unique tracial state was obtained in
\cite{GL2}.

\begin{Cor}\label{MT1C}
Let $X$ be a compact metric space with torsion free $K_i(C(X))$ ($i=0,1$)
and let $A$ be a unital separable simple \CA\, with $TR(A)=0.$
Suppose that $h_1, h_2: C(X)\to A$ are two unital monomorphisms.
Then $h_1$ and $h_2$ are approximately unitarily
equivalent if and only if
$$
(h_1)_{*i}=(h_2)_{*i}\,\,\,\,\,i=0,1
$$
and $\tau\circ h_1(a)=\tau\circ h_2(a)$ for all
$a\in C(X).$
\end{Cor}







We will give two interesting applications of these results. Some
approximate versions of it will be give in Section 4.

\subsection{Minimal dynamical systems}
Let $X$ be a compact metric space and $\alpha: X\to X$ be a
homeomorphism. Recall that $\alpha$ is said to be minimal if
$\{\alpha^n(x):n\in \Z\}$ is dense in $X.$ We assume that $X$ is
an infinite set. The corresponding transformation group \CA\,
denoted by $A_{\alpha}=C(X)\rtimes_{\alpha} {\mathbb Z}$ is a
unital simple \CA. It is also amenable and separable. It is in the
so-called ``Bootstrap" class of \CA s. Therefore it satisfies the
Universal Coefficient Theorem. Many \CA s described above have
tracial rank zero. For example, if $X$ is a connected manifold and
$\alpha$ is a diffeomorphism. Then $TR(A_{\alpha})=0$ if and only
if $\rho_{A_{\alpha}}(K_0(A_{\alpha}))$ is dense in
$Aff(T(A_{\alpha}))$ (see \cite{LqP2} and \cite{Lncre}).

{\it In what follows, we will use $A_{\alpha}$ for
$C(X)\rtimes_{\alpha}\Z$ and $j_{\alpha}: C(X)\to A_{\alpha}$ for
the obvious embedding. }

A theorem of Tomiyama (see \cite{T}) establishes an important
relation between \CA\, theory and topological dynamics:

\begin{thm}{\rm ( J. Tomiyama)} \label{TJT}
Let $X$ be a compact metric space and $\alpha, \beta: X\to X$ be
homeomorphisms. Suppose that $(X, \alpha)$ and $(X, \beta)$ are
topologically transitive. Then $\alpha$ and $\beta$ are flip
conjugate if and only if there is an isomorphism $\phi:
C(X)\rtimes_{\alpha}{\mathbb Z} \to C(X)\rtimes_{\beta}{\mathbb
Z}$ such that $\phi\circ j_{\alpha}=j_{\beta}\circ \chi$ for some
isomorphism $\chi: C(X)\to C(X).$
\end{thm}

It should be noted that all minimal dynamical systems are
transitive.

In the light of Tomiyama's theorem, we introduce the following
version of approximate flip conjugacy (for the case that
$TR(A_{\alpha})=TR(A_{\beta})=0$).

\begin{Def}\label{IIDM}
{\rm Let $(X, \alpha)$ and $(X,\beta)$ be two minimal dynamical
systems such that $TR(A_{\alpha})=TR(A_{\beta})=0.$
 We say that $(X, \alpha)$ and $(X, \beta)$ are {\it
$C^*$-strongly approximately flip conjugate} if there exist
sequences of isomorphisms $\phi_n: A_{\alpha} \to A_{\beta},$
$\psi_n: A_{\beta}\to A_{\beta}$ and sequences of isomorphisms
$\chi_n,\lambda_n: C(X)\to C(X)$ such that
$[\phi_n]=[\phi_1]=[\psi_n^{-1}]$ in $KL(A_{\alpha}, A_{\beta})$
for all $n$ and
$$
\lim_{n\to\infty}\|\phi_n\circ j_{\alpha}(f)-j_{\beta}\circ
\chi_n(f)\|=0 \andeqn \lim_{n\to\infty}\|\psi_n\circ
j_{\beta}(f)-j_{\alpha}\circ \lambda_n(f)\|=0
$$
for all $ f\in C(X).$ }
\end{Def}

For the general case that the crossed products are not assume to
have tracial rank zero, a modified definition is given in
\cite{LM3}. A simplification of teh definition is given below (\ref{IIPDM1}).

In Theorem \ref{TJT}, let $\theta=[\phi]$ in $KK(A_{\alpha},
A_{\beta}).$ Let $\Gamma(\theta)$ be the induced element in\\
$Hom(K_*(A_{\alpha}), K_*( A_{\beta}))$ which preserves the order
and the unit. Then one has
$$
[j_{\alpha}]\times \theta=[j_{\beta}\circ \chi]
$$

Suppose that $TR(A_{\alpha})=TR(A_{\beta})=0.$ Then
$\rho_{A_{\alpha}}(K_0(A_{\alpha}))$ and
$\rho_{A_{\beta}}(K_0(A_{\beta}))$ are dense in
$Aff(T(A_{\alpha}))$ and $Aff(T(A_{\beta})),$ respectively. Thus
$\Gamma(\theta)$ induces an order and unit  preserving affine
isomorphism $\theta_{\rho}: Aff(T(A_{\alpha}))\to
Aff(T(A_{\beta})).$ Recall
$\rho_{A_{\alpha}}:(A_{\alpha})_{s.a}\to Aff(T(A_{\alpha}))$ is
defined by $\rho_{A_{\alpha}}(a)(\tau)=\tau(a)$ for $\tau\in
T(A).$
Therefore, in terms of $K$-theory and $KK$-theory, one has the
following: If $\alpha$ and $\beta$ are flip conjugate, then there
is an isomorphism $\chi:C(X)\to C(X)$ such that
\begin{eqnarray}\label{eEK1}
[j_{\alpha}]\times \theta=[j_{\beta}\circ \chi]\,\,\,{\rm  in}
\,\,\, KK(C(X), A_{\beta})\,\,{\rm and}\,\,\theta_{\rho}\circ
\rho_{A_{\alpha}}\circ j_{\alpha}= \rho_{A_{\beta}}\circ
j_{\beta}\circ \chi\,\,\,{\rm on}\,\,\,C(X)_{s.a}.
\end{eqnarray}

The following theorem gives a $K$-theoretical description of
$C^*$-strong approximate flip conjugacy:

\begin{thm}\label{IIITM1}
Let $(X, \alpha)$ and $(X,\beta)$ be two minimal dynamical systems
such that $A_{\alpha}$ and $A_{\beta}$ have tracial rank zero.
Then $\alpha$ and $\beta$ are $C^*$-strongly approximately flip
conjugate if and only if the following hold: There are  sequences
of isomorphisms  $\chi_n,\lambda_n: C(X)\to C( X)$ and $\theta\in
KL(A_{\alpha}, A_{\beta})$ such that  $\Gamma(\theta)$ gives an
isomorphism from $(K_0(A_{\alpha}),K_0(A_{\alpha})_+, [1],
K_1(A_{\alpha}))$ to $(K_0(A_{\beta}),K_0(A_{\beta})_+, [1],
K_1(A_{\beta})),$ for any finitely generated subgroup $G\subset
\underline{K}(C(X)),$
\begin{eqnarray}\label{eDy1}
[j_{\alpha}]\times \theta|_G=[j_{\beta}\circ \chi_n]|_G \andeqn
[j_{\beta}]\times \theta^{-1}|_G=[j_{\alpha}\circ \lambda_n]|_G,
\end{eqnarray}
for all sufficiently large $n$ and
\begin{eqnarray}\label{eDy2}
 \lim_{n\to\infty}\|\rho_{A_{\beta}}\circ j_{\beta}\circ
\chi_n(f)-\theta_{\rho}\circ \rho_{A_{\alpha}}\circ
j_{\alpha}(f)\|=0\andeqn
\end{eqnarray}
\begin{eqnarray}\label{eDye2+}
\lim_{n\to\infty}\|\rho_{A_{\alpha}}\circ j_{\alpha}\circ
\lambda_n(f)-\theta_{\rho}^{-1}\circ \rho_{A_{\beta}}\circ
j_{\beta}(f)\|=0
\end{eqnarray}
for all $ f\in C(X)_{s.a}.$

\end{thm}

\begin{Cor}\label{IICapp}
Let $X$ be a compact metric space with torsion free $K$-theory.
Let $(X, \alpha)$ and $(X, \beta)$ be two minimal dynamical
systems such that $TR(A_{\alpha})=TR(A_{\beta})=0.$ Suppose that
there is a unit preserving order isomorphism
\begin{eqnarray}\label{eDy3}
\gamma: (K_0(A_{\alpha}), K_0(A_{\alpha})_+, [1_{A_{\alpha}}],
K_1(A_{\alpha})) \to (K_0(A_{\beta}), K_0(A_{\beta})_+,
[1_{A_{\beta}}], K_1(A_{\beta})),
\end{eqnarray}
\begin{eqnarray}\label{eDy4}
[j_{\alpha}]\times \theta=[j_{\beta}\circ \chi] \,\,\,{\rm
in}\,\,\, KL(C(X), A_{\beta})\,\,{\rm and}\,\, \gamma_{\rho}\circ
j_{\alpha}=\rho_{A_{\beta}}\circ j_{\beta}\circ \chi\,\,{\rm
on}\,\, C(X)_{s.a}
\end{eqnarray}
for some isomorphism $\chi: C(X)\to C(X).$ Then $(X, \alpha)$ and
$(X, \beta)$ are $C^*$-strongly approximately flip conjugate.
\end{Cor}

\begin{Remark}\label{Rcoj}
{\rm In definition \ref{IIDM}, we require that
$[\phi_n]=[\phi_1]$ for all $n.$ It will  be only a marginal gain
by removing  this condition from the definition. If both
$K_i(A_{\alpha})$ and $K_i(A_{\beta})$ are finitely generated
($i=0,1$), then there are only finitely many order isomorphisms
which preserve the identity. Moreover, in this case,
$ext_{\Z}(K_{i-1}(A_{\alpha}),K_i(A_{\beta}))$ has only finitely
many elements. Thus, in this case, there are only finitely many
elements $\theta\in KL(A_{\alpha}, A_{\beta})$ which gives order
and unit preserving isomorphisms from $K_i(A_{\alpha})$ to
$K_i(A_{\beta})$ ($i=0,1$). Hence, there is a subsequence
$\{n_k\}$ such that $[\phi_{n_k}]=\theta$ for some $\theta.$
Therefore  one may well assume that $[\phi_n]=[\phi_1]$ for all
$n.$

In the case when $X$ is the Cantor set, $K_0(C(X))=C(X, \Z).$ It
follows that conditions in \ref{IIITM1} can be considerably
simplified. 
Furthermore, if $X$ is the Cantor set, two minimal homeomorphisms
$\alpha$ and $\beta$ are $C^*$-strongly approximately conjugate if
and only if they are approximately $K$-conjugate (see \ref{IIIDconjK}). More precisely
we have the following theorem which is also related to the work of
Giordano, Putnam and  Skau in  \cite{GPS}.
}
\end{Remark}
\begin{thm}{\rm ( Theorem 5.4 of \cite{LM})}\label{TCant}
Let $X$ be the Cantor set and $\alpha$ and $\beta$ be
minimal homeomorphisms.
Then the following are equivalent:

{\rm (i)} $\alpha$ and $\beta$ are $C^*$-strongly approximately
flip conjugate,

{\rm (ii)} $\alpha$ and $\beta$ are approximately $K$-conjugate,

{\rm (iii)} $A_{\alpha}$ and $A_{\beta}$ are isomorphic,

{\rm (iv)} $(K_0(A_{\alpha}), K_0(A_{\alpha})_+, [1_{A_{\alpha}}])
\cong (K_0(A_{\beta}), K_0(A_{\beta})_+, [1_{A_{\beta}}]),$

{\rm (v)} there exists a sequence $\gamma_n\in [[\alpha]]$ and
$\sigma_n\in [[\beta]]$ and a homeomorphism $\chi: X\to X$ such
that
$$
f\circ \alpha=\lim_{n\to\infty}
f\circ \chi\circ \sigma_n\circ \beta\circ \sigma_n^{-1}\circ \chi^{-1}
\Tand
$$
$$
f\circ\beta=\lim_{n\to\infty}f\circ \chi^{-1}\circ \gamma_n\circ \alpha \circ
\gamma_n^{-1}\circ \chi
$$
for all $f\in C(X).$

{\rm (vi)} $\alpha$ and $\beta$ are strong orbit equivalent.

\end{thm}

In Theorem \ref{TCant}, $[[\alpha]]$ is the set of topological full
group with respect to $\alpha,$ i.e., the group of all
homeomorphisms $\gamma: X\to X$ such that
$\gamma(x)=\alpha^{n(x)}(x)$ for all $x\in X,$ where $n\in C(X,
\Z).$ We will explain other terminologies used in  Theorem
\ref{TCant} in section 5. Further discussion of approximate
conjugacy will be given in section 5.

\subsection{$C^*$-dynamical systems and the Rokhlin property}

By a $C^*$-dynamical system we mean a pair $(A,\alpha),$ where $A$
is a \CA\, and $\alpha\in Aut(A).$ Let $A$ be a unital simple
$A\T$-algebra with real rank zero and $\alpha\in Aut(A)$ be a
``sufficiently outer" automorphism. A. Kishimoto studied the
problem when the associated  crossed product is again an
$A\T$-algebra of real rank zero (\cite{Ks1},
\cite{Ks1.5},\cite{Ks2}  and \cite{Ks3}). In particular, Kishimoto
studied the case that $A$ has a unique tracial state and $\alpha$
is approximately inner. Kishimoto also suggested that the
appropriate notion for ``sufficiently outer" is the Rokhlin
property. A more general question is: Let $A$ be a unital simple
AH-algebra with no dimension growth and with real rank zero.
Suppose that $\alpha\in Aut(A).$ When is $A\rtimes_{\alpha} \Z$
again a unital simple AH-algebra with no dimension growth and with
real rank zero?

With the classification theorem for simple
\CA s of tracial topological rank zero, an important question
related to the crossed products is the following: Let $A$ be a
unital separable simple \CA\, with tracial rank zero and $\alpha$
is an (outer) automorphism. When does $A\rtimes_{\alpha}\Z$ have
tracial rank zero?

The following is defined in \cite[Definition 2.1]{OP}.

\begin{Def}\label{DRok1}
{\rm
Let $A$ be a simple unital \CA~ and
let $\alpha \in \Aut(A)$.
We say that $\alpha$ has the {\it tracial Rokhlin property}
if for every finite set $F \subset A$, every $\ep > 0$,
every $n \in {\mathbb N}$, and every nonzero positive element $x \in A$,
there are mutually orthogonal projections $e_0, e_1, \dots, e_n \in A$
such that:
\begin{itemize}
\item[$(1)$]
$\| \alpha (e_j) - e_{j + 1} \| < \ep$ for $0 \leq j \leq n - 1$.
\item[$(2)$]
$\| e_j a - a e_j \| < \ep$ for $0 \leq j \leq n$ and all $a \in F$.
\item[$(3)$]
With $e = \sum_{j = 0}^{n} e_j$, $[1 - e] \leq [x]$.
\end{itemize}
}

\end{Def}

Recall that a unital simple \CA\, $A$ is said to have the
Fundamental Comparison Property, if for any two projections $p,
q\in A$ $p\sim q'\le q,$ if $\tau(p)<\tau(q)$ for all $\tau\in
T(A).$ In the definition \ref{DRok1}, if $A$ has the Fundamental
Comparison Property, then condition (3) can be replaced by
$$
\hspace{-2.68in}{ \rm (3)'}\,\,\,\,\,\,\,\,\,\,\,\,\hspace{0.6in}
\tau(1-e)<\ep\rforal \tau\in A.
$$
This is the reason why it called tracial Rokhlin property. A much
stronger version of Rokhlin property was introduced in \cite{LO}
recently.

\begin{Def}\label{DRokc}
{\rm Let $A$ be a simple unital
 \CA~ and
let $\alpha \in \Aut(A)$.
We say $\alpha$ has the {\it tracial cyclic Rokhlin property}
if for every finite set $F \subset A$, every $\ep > 0$,
every $n \in {\mathbb N}$, and every nonzero positive element $x \in A$,
there are mutually orthogonal projections $e_0, e_1, \dots, e_n \in A$
such that
\begin{itemize}
\item[$(1)$]
$\| \alpha (e_j) - e_{j + 1} \| < \ep$ for $0 \leq j \leq n$,
where $e_{n + 1} = e_0$.
\item[$(2)$]
$\| e_j a - a e_j \| < \ep$ for $0 \leq j \leq n$ and all $a \in F$.
\item[$(3)$]
With $e = \sum_{j = 0}^{n} e_j$, $[1 - e] \leq [x]$.
\end{itemize}
}
\end{Def}

Note that if $\alpha$ has the tracial Rokhlin property, then
$K_0(A)$ must have a ``dense" subset which is invariant under
$\alpha_{*0}.$ It is shown (see Theorem 2.9 of \cite{LO}) that if
$\alpha$ has tracial cyclic Rokhlin property then indeed the
crossed product $A\rtimes_{\alpha}\Z$ has tracial rank zero. It is
proved in \cite{OP} that tracial Rokhlin property occurs more
often than one may think. For example, assuming  that $A$ is a
unital separable simple \CA\, with a unique tracial state and with
tracial rank zero, they prove in \cite{OP} that
$A\rtimes_{\alpha}\Z$ has tracial Rokhlin property if and only if
$A\rtimes_{\alpha}\Z$ has a unique tracial state, or
$A\rtimes_{\alpha}\Z$ has real rank zero.
 It is
proved in \cite{LO} that, if $\alpha^r$ is approximately inner for
some integer $r>0$ and $\alpha$ has the tracial Rokhlin property,
then $\alpha$ has the tracial cyclic Rokhlin property.
Consequently, $A\rtimes_{\alpha}\Z$ has tracial rank zero. In
particular, by the classification theorem for simple unital
separable amenable \CA s with tracial rank zero (see
\cite{Lnduke}), if  $A$ is a unital $A\T$-algebra of real rank
zero, in this case, $A\rtimes_{\alpha}\Z$ is again a unital
$A\T$-algebra with real rank zero provided that
$A\rtimes_{\alpha}\Z$ has torsion free $K$-theory. This solves the
Kishimoto's problem in this setting. By applying an approximate
version of theorem \ref{MT02} (see \ref{MT1C1}), we have the
following theorem:

\begin{thm}\label{MTrok}
Let $A$ be a unital separable simple amenable \CA\, with $TR(A)=0$
which satisfies the UCT and let $\alpha\in Aut(A).$ Suppose that
$\alpha$ satisfies the tracial Rokhlin property. If there is an
integer $r>0$ such that $\alpha^r_{*0}|_G={\rm id}_G$ for some
subgroup $G\subset K_0(A)$ for which $\rho_A(G)=\rho_A(K_0(A)),$
then $\alpha$ satisfies the tracial cyclic Rokhlin property.
\end{thm}

From this we obtain:

\begin{thm}\label{MTAH}
Let $A$ be a unital separable simple amenable \CA\, with $TR(A)=0$
which satisfies the UCT and let $\alpha\in Aut(A).$ Suppose that
$\alpha$ satisfies the tracial Rokhlin property and
$\alpha^r_{*0}|_G={\rm id}_{G}$ for some subgroup $G\subset
K_0(A)$ for which $\rho_A(G)=\rho_A(K_0(A)).$ Then
$A\rtimes_{\alpha}\Z$ is a unital simple AH-algebra with no
dimension growth and with real rank zero.
\end{thm}

As a consequence, combining the results in \cite{OP},  we also
have the following:

\begin{thm}\label{VTRC}
Let $A$ be a unital separable simple amenable \CA\, with $TR(A)=0$
and with a unique tracial state which satisfies the UCT and let
$\alpha\in Aut(A)$ be such that $\alpha^k_{*0}|G={\rm id}_{G}$ for
some subgroup $G\subset K_0(A)$ for which $\rho_A(G)=\rho_A(
K_0(A)).$ Then the following are equivalent:

{\rm (i)} $\alpha$ has the tracial Rokhlin property;

({\rm ii)} $A\rtimes_{\alpha}{\mathbb Z}$ has a unique tracial
state;

{\rm (iii)} $A\rtimes_{\alpha}{\mathbb Z}$ has real rank zero;

{\rm (iv)} $\alpha$ has  the tracial cyclic Rokhlin property;

{\rm (v)} $A\rtimes_{\alpha}\Z$ has tracial rank zero;

{\rm (vi)} $A\rtimes_{\alpha}\Z$ is a unital simple AH-algebra
with unique tracial state and real rank zero.

\end{thm}

Finally we would like to mention the Furstenberg transformations
on irrational rotation algebras studied recently by H. Osaka and
N. C. Phillips. These also include the transformation group \CA s
of  minimal Furstenberg transformations on the torus.

\begin{Def}\label{DFur}
{\rm
Let $A_{\theta}$ be the usual rotation algebra generated by
unitaries $u$ and $v$ satisfying $vu=e^{2\pi i\theta}uv.$ Let
$\theta, \gamma\in \R,$ and let $d\in \Z,$ and let $f: S^1\to \R$
be a continuous function. The {\it Furstenberg transformation} on
$A_{\theta}$ determined by $(\theta, \gamma, d, f)$ is the
automorphism $\alpha_{\theta, \gamma, d,f}$ of $A_{\theta}$ such
that
$$
\alpha_{\theta, \gamma, d, f}(u)=e^{2\pi i\gamma} u\andeqn
\alpha_{\theta,\gamma, d,f}(v)=\exp(2\pi i f(u))u^dv.
$$
The resulted crossed product is denoted by
$A_{\theta}\rtimes_{\alpha}\Z.$ }
\end{Def}

Combining the results in \cite{OP} and Theorem \ref{VTRC}, we
obtain the following.

\begin{thm}\label{TLFur}
Let $\theta, \gamma\in \R$ and suppose that $1, \theta,\gamma$ are
linearly independent over $\Q.$ Let $d\in \Z$ and let
$\alpha_{\theta,\gamma,d,f}\in Aut(A_{\theta})$ be as defined in
\ref{DFur}. Then

{\rm (1)} $A_{\theta}\rtimes_{\alpha}\Z$ is a unital simple
AH-algebra with no dimension growth,  with real rank zero and with
a unique tracial state.

{\rm (2)} If, in addition, $d=0,$ then
$A_{\theta}\rtimes_{\alpha}\Z$ is an A$\T$-algebra with real rank
zero.
\end{thm}

\section{Maps from $C(X)$}

In this section we will prove  Theorem \ref{MT02} and
Theorem \ref{MT1}. We do this by proving an approximate version of
these (Theorem \ref{MT01}).

\begin{Lem}\label{Lapt}
Let $A$ be a unital  \CA.
For any $\ep>0$ and any finite subset ${\cal F}\subset A,$
there exists $\dt>0$ and ${\cal G}\subset A$ which satisfy
the following:

If $B$ is another unital \CA,  $\tau\in T(B)$ and
$\phi_n: A\to B$ is a unital \morp\, which is
${\cal G}$-$\dt$-multiplicative, then
there exists a tracial state $\sigma\in T(A)$
such that
$$
|\tau\circ \phi_n(a)-\sigma(a)|<\ep \rforal a\in {\cal F}.
$$

\end{Lem}

\begin{proof}
Suppose that the lemma is false. Then there exist $\ep_0>0$ and a
finite subset ${\cal F}_0\subset A,$ there exists a ${\cal
G}_n$-$\dt_n$-multiplicative \morp\, $\phi_n: A\to B_n,$ where
$B_n$ are unital \CA s and $\cup_{n=1}^{\infty}G_n$ is dense in
$A$ and $\sum_{n=1}^{\infty}\dt_n<\infty,$ and a sequence of
tracial states $\tau_n\in T(B_n)$ such that
$$
\inf_{n\in N}\{\inf\{\sup\{|\tau_n\circ \phi_n(a)-\sigma(a)|: a\in
{\cal F}_0\}: \sigma\in T(A)\}\}\ge \ep_0.
$$
Take a weak limit $\sigma$ of $\{\tau_n \circ \phi_n\}.$ Then
$\sigma$ is a tracial state of $A.$
There is subsequence $\phi_{n(k)}$ such
that
$$
\sigma(a)=\lim_{k\to\infty}\tau\circ \phi_{n(k)}(a)\rforal a\in A.
$$
This gives a contradiction.
\end{proof}

\begin{Def}\label{Dnota}
{\rm Let $X$ be a compact metric space and let $A$ be a unital \CA. Suppose
that $\phi: C(X)\to A$ is a unital positive linear map and suppose
that $\tau\in T(A).$ Then $\tau\circ \phi$ is a positive linear
functional on $C(X).$ We use  $\mu_{\tau\circ \phi}$
for the induced probability Borel regular measure.

Let $X$ be compact metric space and $\eta>0.$
Then there are finitely many (distinct) points $x_1,x_2,...,x_m\in X$ such that
$\{x_1,x_2,...,x_m\}$ is an $\eta$-dense subset.
There is an integer $s>0,$ such that
$$
O_i\cap O_j=\emptyset,\,\,\,\,\,i\not=j,
$$
where $O_i=\{x\in X: {\rm dist}(x, x_i)<\eta/2s\},$ $i=1,2,...,m.$
The integer $s$ depends on the choice of $\{x_1,x_2,...,x_m\}.$
}
\end{Def}

\begin{Lem}\label{Ldig2}
Let $X$ be a compact metric space, $\ep>0$
and ${\cal F}\subset
C(X)$ be a finite subset. Let $L>0$ be an integer
and let $\eta>0$ be such that
$|f(x)-f(x')|<\ep/8$ if ${\rm dist}(x,x')<\eta.$
Then, for any integer $s>0,$
any finite $\eta/2$-dense subset $\{x_1,x_2,...,x_m\}$ of $X$
for which $O_i\cap
O_j=\emptyset,\,\,{\rm if}\,\,\, i\not=j,
$ where
$$
O_i=\{\xi\in X: {\rm dist}(\xi, x_i)<\eta/2s\}
$$
and any  $1/2s>\sigma>0,$ there exist
a finite subset ${\cal G}\subset C(X)$ and $\dt>0$
satisfying the following:

For any unital separable stably finite \CA\, $A$ with real rank
zero, $\tau\in T(A)$ and any ${\cal G}$-$\dt$-multiplicative
\morp\, $\phi: C(X)\to A,$
if $\mu_{\tau\circ \phi}(O_i)>\sigma\cdot
\eta,$ for all $i,$
then there are mutually orthogonal
projections $p_1,p_2,...,p_m$ in $A$ such that
$$
\|\phi(f)-((1-p)\phi(f)(1-p)+\sum_{i=1}^mf(x_i)p_i)\|<\ep\rforal
f\in {\cal F}
$$
$$
 \andeqn \|(1-p)\phi(f)-\phi(f)(1-p)\|<\ep,
$$
where $p=\sum_{i=1}^mp_i$
$$
\tau(p_k)>(4mL+4)\tau(1-\sum_{i=1}^m p_i)\andeqn
\tau(p_k)>{\sigma\over{2}}\cdot \eta,\,\,k=1,2,...,m.
$$
\end{Lem}


\begin{proof}
Suppose
that there exists $\ep_0>0$ and a finite subset ${\cal F}_0\subset C(X)$
so that the lemma is false.  Fix $\ep>0$ so that $\ep<\ep_0/4.$
Let $\eta>0$ so that
$$
|f(x)-f(y)|<\ep/8\,\,\,{\rm if}\,\,\, {\rm dist}(x, y)<\eta,\,\,\, x,
y\in X
$$
for all $f\in {\cal F}_0.$ Suppose that $\{x_1,x_2,...,x_m\}$ is
an $\eta/2$-dense subset of $X$ such that $O_i\cap O_j=\emptyset,$
if $i\not=j,$ where
$$
O_i=\{x\in X: {\rm dist}(x,x_i)<\eta/2s\},\,i=1,2,...,m
$$
for some integer $s>0.$
Let $1/2s>\sigma>0.$
Since the lemma is
false (for a choice of
the aboved mentioned $x_1,...,x_m,$ $s$ and $\sigma$),
there  exist a sequence of unital separable stably finite
\CA\, $B_n$ of real rank zero, a sequence of unital ${\cal
G}_n$-$\dt_n$-multiplicative \morp s $\phi_n: C(X)\to B_n,$
where $\cup_{n=1}^{\infty} G_n$ is dense in $C(X)$ and
$\sum_{n=1}^{\infty}\dt_n<\infty,$
and there are  $\tau_n\in T(B_n)$ such that $\mu_{\tau\circ
\phi_n}(O_i)>\sigma\cdot \eta$ satisfying the following:
$$
\inf\{\sup\{\|\phi_n(f)-[(1-p_n)\phi_n(1-p_n)+\sum_{i=1}^{m}f(x_i)p(i,n)]\|
: f\in {\cal F}_0\}: n\}\ge \ep_0
$$
where the infimum is taken among all mutually orthogonal
projections $p(1,n),p(2,n),..., p(m,n)$ which satisfy
$$
\tau_n(p(i,n))>(4mL+4)\tau_n(1-p_n) \andeqn \tau_n(p(i,n))>\sigma\eta/2,
$$
where $p_n=\sum_{i=1}^{m}p(i,n).$

Define $\Phi: C(X)\to l^{\infty}(\{B_n\})$ by
$\Phi(f)=\{\phi_n(f)\}$ and let $\pi: l^{\infty}(\{B_n\})\to
q_{\infty}(\{B_n\}).$ Then $\pi\circ \Phi: C(X)\to
q_{\infty}(\{B_n\})$ is a \hm. By passing to a subsequence, if
necessary, we may also assume that there is $\tau\in
l^{\infty}(\{B_n\})$ such that
$$
\lim_{n\to\infty}\tau_n(a_n)=\tau(\{a_n\})
$$
for  any $\{a_n\}\in l^{\infty}(\{B_n\}).$
Moreover, since, for each
$\{a_n\}\in c_0(\{B_n\}),$
$\lim_{n\to\infty}\tau_n(a_n)=0,$ we may view $\tau$ as a tracial
state of $q_{\infty}(\{B_n\}).$
It follows that
$$
\mu_{\tau,\pi\circ\Phi}(O_i)>\sigma\cdot \eta,\,\,\,\,\,
i=1,2,...,m.
$$
Exactly as in the proof of 2.11 of \cite{Lnfields}, there are
$r_i\in (0,1)$ such that
\begin{eqnarray}\label{eC}
\mu_{\tau,\pi\circ \Phi}(C_{r_i})=0,
\end{eqnarray}
where
$$
C_{r_i}=\{\zeta\in X: {\rm dist}(\zeta, x_i)= (1+r_i)(\eta/2)\}.
$$
Let $\Omega_1, \Omega_2,...,\Omega_{K'}$ be
disjoint open subsets such that $\cup_{i=1}^{K'}\Omega_i=X\setminus
\cup_{i=1}^mC_{r_i}.$ Therefore
$$
\cup_{i=1}^{K'}\Omega_i\cup (\cup_{i=1}^mC_{r_i})=X\andeqn {\rm
diam}(\Omega_i)<\eta/2.
$$
(Note  that $K'\le m^m$.)
Since $O_i\cap O_j=\emptyset,$ if $i\not=j,$
 we may assume that $O_i\subset \Omega_i,$
$i=1,2,...,m$ and $\mu_{\tau\circ \pi\circ \Phi}(\Omega_i)>0,$
$i=1,2,...,K$ (and $\mu_{\tau\circ \pi\circ \Phi}(\Omega_i)=0$ if $i>K$
and $K'\ge K\ge m$).
Let $B_{\Omega_i}={\rm Her}(\pi\circ \Phi(\Omega_j)),$ $i=1,2,...,K'.$
Since $B_n$ has real rank zero,
$q_{\infty}(\{B_n\})$ also has real rank zero. Hence there is an
approximate identity $\{e_i(n)\}$ for $B_{O_i}$ for
each $i.$
Then
$$
\tau(e_i^{(n)})\nearrow \mu_{\tau\circ \pi\circ \Phi}(\Omega_i),
$$
$i=1,2,...,K.$
Since $\mu_{\tau\circ\pi \circ \Phi}(X\setminus \cup_i^m C_{r_i})=0,$
$$
\tau(\sum_{i=1}^Ke_i^{(n)})\nearrow 1
$$
as $n\to\infty.$ So
$$
\tau(1-\sum_{i=1}^K e_i^{(n)})\to 0\,\,\,{\rm as}\,\,\,n\to\infty.
$$
Since $\mu_{\tau\circ \pi\circ \Phi}(\Omega_i)>0$ for $i=1,2,...,K,$
we obtain projections $q_j$ ($=e_j(n)$ for some
large $n$) so that
$$
(4KL+5)\tau(1-\sum_{k=1}^Kq_k)<\tau(q_i)\,\,\,i=1,2,...,K.
$$
By Lemma 2.5  of \cite{Lnfields} and by the choice of $\eta,$ one
obtains
$$
\|\pi\circ \Phi(f)-[(1-q)\pi\circ
\Phi(f)(1-q)+\sum_{k=1}^Kf(\zeta_k)q_k]\|<\ep/2 \andeqn
$$
$$
\|q\pi\circ \Phi(f)-\pi\circ \Phi(f)q\|<\ep/2
$$
for all $f\in {\cal F}_0,$ where $q=\sum_{k=1}^Kq_k.$
Thus, for all sufficiently large $n,$ there are nonzero mutually
orthogonal projections $p(1,n),p(2,n),...,p(K,n)\in B_n$ such
that
$$
\|\phi_n(f)-[(1-p_n)\phi_n(f)(1-p_n)+\sum_{k=1}^Kf(\zeta_k)p(k,n)]\|<\ep/2,
$$
$$
\|p_n\phi_n(f)-\phi_n(f)p_n\|<\ep/2\rforal f\in {\cal F}_0 \andeqn
$$
$$
(4KL+4)\tau_n(1-p_n)<\tau_n(p(k,n)),
$$
where $p_n=\sum_{k=1}^Kp(i,n).$ Note that $O_i\subset \Omega_i,$
$i=1,2,...,m.$ Since $\{x_1,x_2,...,x_m\}$ is $\eta/2$-dense in
$X,$ by the choice of $\eta$ and by replacing some of points
$\xi_j$ which close to $x_i$ within $\eta/2$ by $x_i$ and $p(k,n)$
by sum of some $p(k',n)$'s,
 we may assume that $K=m$ and $\xi_i=x_i.$ Furthermore, it
is also clear that we may assume that
$$
\tau(p_k)>{\sigma\cdot\eta\over{2}}, \,\,\,k=1,2,...,m.
$$
 This is a contradiction.
\end{proof}

\begin{Lem}\label{L2dig}
Let $X$ be a compact metric space, $\ep>0$ and ${\cal F}\subset
C(X)$ be a finite subset. Let $L>0$ be an integer and
$\eta>0$ be such that $|f(x)-f(x')|<\ep/8$ if ${\rm dist}(x,x')<\eta.$
 Then,
for any integer $s>0,$ any finite $\eta/2$-dense subset of $X$ for which
$O_i\cap O_j=\emptyset$ for $i\not=j,$ where
$$
O_i=\{\xi\in X: {\rm dist}(\xi, x_i)<\eta/2s\}
$$
and any $1/2s>\sigma>0,$
there exist $\gamma>0,$
a finite subset ${\cal G}\subset C(X)$ and $\dt>0$ satisfying the
following:

For any unital separable stably finite \CA\, $A$ with real rank
zero, $\tau\in T(A)$ and any unital ${\cal G}$-$\dt$-multiplicative
\morp\, $\phi, \psi: C(X)\to A$ with
$$
|\tau\circ \phi(g)-\tau\circ \psi(g)|<\gamma\,\,\, \rforal g\in
{\cal G},
$$
 and if
$$
\mu_{\tau\circ \phi}(O_i),\,\,\, \mu_{\tau\circ \psi}(O_i)
>{\sigma\cdot \eta},
$$
for all $i,$  then there are mutually
orthogonal projections $p_1,p_2,...,p_m$ and $q_1,q_2,...,q_m$ in
$A$ such that
$$
\|\phi(f)-((1-p)\phi(f)(1-p)+\sum_{i=1}^mf(x_i)p_i)\|<\ep\rforal
f\in {\cal F}
$$
\vspace{-0.05in}
$$
 \andeqn \|(1-p)\phi(f)-\phi(f)(1-p)\|<\ep/2,
$$
\vspace{-0.05in}
 where $p=\sum_{i=1}^mp_i,$
$$
\tau(p_k)>(3mL+3)\tau(1-\sum_{i=1}^m p_i)\andeqn
\tau(p_k)>{\sigma\cdot \eta\over{2}},k=1,2,...,m,\,\,\,\,\,\andeqn
$$
\vspace{-0.05in}
$$
\|\psi(f)-((1-q)\psi(f)(1-q)+\sum_{i=1}^mf(x_i)q_i)\|<\ep\rforal
f\in {\cal F}
$$
\vspace{-0.05in}
 $$
 \andeqn \|(1-q)\psi(f)-\psi(f)(1-q)\|<\ep,
$$
where $q=\sum_{i=1}^mq_i,$
$$
\tau(q_k)>(3mL+3)\tau(1-\sum_{i=1}^m q_i)\andeqn
\tau(q_k)>{\sigma\cdot \eta/2},\,\,\,k=1,2,...,m,\andeqn
$$
\vspace{-0.05in}
\begin{eqnarray}\label{ede} |\tau(q_k)-\tau(p_k)|<
{\min_k\{\tau(p_k)\}\over{3mL}}, \,\,\,k=1,2,...,m.
\end{eqnarray}

\end{Lem}


\begin{proof}
First we note that only (\ref{ede}) needs a proof. The proof is
similar to that of \ref{Ldig2}. But we need to ``dig" out the
projections simultaneously.

Let $B_n$ be any sequence of unital separable \CA s of real rank
zero, $\phi_n, \psi_n: C(X)\to B_n$ be any unital \morp s  and
$\tau_n\in T(B_n)$  such that
$$
\lim_{n\to\infty}\|\phi_n(ab)-\phi_n(a)\phi_n(b)\|=0,
\lim_{n\to\infty}\|\psi_n(ab)-\psi_n(a)\psi_n(b)\|=0\andeqn
$$
$$
\lim_{n\to\infty}|\tau_n\circ \phi_n(f)-\tau_n\circ \psi_n(f)\|=0
$$
for all $a,b, f\in C(X).$ Let $\Phi, \Psi: C(X)\to
l^{\infty}(\{B_n\})$ be defined by $\Phi(f)=\{\phi_n(f)\}$ and
$\Psi_n(f)=\{\psi_n(f)\},$ respectively. Let $\pi:
l^{\infty}(\{B_n\})\to q_{\infty}(\{B_n\})$ be the quotient map.
As in  the proof of \ref{Ldig2}, $\pi\circ \Phi, \pi\circ \Psi:
C(X)\to q_{\infty}(\{B_n\})$ are unital \hm s. By the assumption,
borrowing the notation in the proof of \ref{Ldig2}, $\tau\circ
\pi\circ \Phi=\tau\circ \pi\circ \Psi.$ Thus, in the proof of
\ref{Ldig2}, one has
$$
\mu_{\tau,\pi\circ \Phi}(\Omega_i)=\mu_{\tau, \pi\circ
\Psi}(\Omega_i),\,i=1,2,...,K.
$$
Note that we assume that $\mu_{\tau,\pi\circ \Phi}(\Omega_i)>0,$
$i=1,2,...,K.$ Note also that $m\le K\le m^m.$ Let $\{e'_j(n)\}$
be an approximate identity for $Her(\pi\circ \Psi_n(\Omega_j)).$
Put
$$
r=\inf\{\mu_{\tau,\pi\circ \Phi}(\Omega_i): i=1,2,...,K\}>0.
$$
So one can choose $n$ so that
$$
|\tau(e'_j(n))-\tau(e_j(n))|<r/(4KL+4) \,\,\,{\rm as\,\,\,
well\,\,\, as}\,\,\,
$$
$$
\tau(e_j'(n))>(4KL+4)\tau(1-\sum_{j=1}^Ke_j'(n)), j=1,2,...,K.
$$
We then apply the last argument of the proof of \ref{Ldig2}.  It
is then clear from the proof of \ref{Ldig2} that, by matching the
size of projections that one may further require  (\ref{ede}) to be
held.
\end{proof}

The following is taken from Theorem 3.1 in \cite{GL3} (see also
Remark 1.1 in \cite{GL3}). Some special cases of this can be found
in \cite{Lnrr0}, \cite{EGLP}. A different form, but similarly in nature,
of the following  was proved in
\cite {D}.

\begin{thm}\label{Luniq}
Let $X$ be a compact metric space. For any $\ep>0$ and any finite
subset ${\cal F}\subset C(X),$ there exist $\dt>0,$ $\eta>0,$ an
integer $N>0,$ a finite subset ${\cal G}\subset C(X)$ and a finite
subset ${\cal P}\subset {\bf P}(C(X))$ satisfy the following:

For any unital \CA\, $A$ with real rank zero, stable rank one and
weakly unperforated $K_0(A)$ and any unital  ${\cal
G}$-$\dt$-multiplicative \morp s $\phi, \psi: C(X)\to A,$ if
$$
[\phi]|_{\cal P}=[\psi]|_{\cal P},
$$
then there exists a unitary $u\in M_{Nk+1}(A)$ such that
$$
\phi(f)\oplus f(x_1)\cdot 1_N\oplus f(x_2)\cdot 1_N\oplus\cdots
f(x_k)\cdot 1_N\aueeps {u}_{\ep}\psi(f), f(x_1)\cdot \oplus 1_N\oplus
f(x_2)\cdot 1_N\oplus\cdots \oplus f(x_k)\cdot 1_N
$$
for all $f\in {\cal F},$ and for any $\eta$-dense set
$\{x_1,x_2,...,x_k\}$ in $X.$
\end{thm}

\begin{thm}\label{MT01}
Let $X$ be a compact metric space, $\ep>0$ and ${\cal F}\subset
C(X)$ be a finite subset. Let $\eta>0$ be
such that $|f(x)-f(x')|<\ep/8$ if ${\rm dist}(x,x')<\eta.$
Then, for any integer $s>0,$ any finite
$\eta/2$-dense subset $\{x_1,x_2,...,x_m\}$ of $X$ for which
$O_i\cap O_j=\emptyset,$ where
$$
O_i=\{x\in X: {\rm dist}(x, x_i)<\eta/2s\}
$$
and any $1/2s>\sigma>0,$  there exist $\gamma>0,$
a finite subset
${\cal G}\subset C(X),$ $\dt>0,$ and a finite subset ${\cal
P}\subset {\bf P}(C(X))$ satisfying the following:

For any unital separable simple \CA\, $A$ with $TR(A)=0,$
any ${\cal G}$-$\dt$-multiplicative unital
\morp s $\phi, \psi: C(X)\to A$ with
$$
|\tau\circ \phi(g)-\tau\circ \psi(g)|<\gamma \rforal g\in {\cal
G} \andeqn \rforal \tau\in T(A), {\rm if}
$$
$$
\mu_{\tau\circ \phi}(O_i),\,\,\, \mu_{\tau\circ \psi}(O_i)
>{\sigma\cdot \eta},
$$
for all $i$ and for all $\tau\in T(A),$
and, if
$$
[\phi]|_{{\cal P}}=[\psi]|_{{\cal P}},
$$
then there exists a unitary $u\in A$ such that
$$
\phi\aueeps{u}_{\ep} \psi\,\,\,\,{\rm on}\,\,\,{\cal F}.
$$

\end{thm}


\begin{proof}
 Fix $\ep>0$ and a finite subset ${\cal F}\subset C(X).$ Let
$\dt_1,$ ${\cal G}_1,$ ${\cal P}_1,$ $\eta_1$ and integer $N$ be
required by \ref{Luniq} for $\ep/8$ and ${\cal F}.$ Let $L=2N.$
There is a finite subset ${\cal F}_1\subset C(X)$ and
$\dt_2>0$ satisfying the following:
for any ${\cal F}_1$-$\dt_2$-multiplicative
\morp s $H_1, H_2: C(X)\to B$ (for any unital $B$)
if $H_1\approx_{\dt_2} H_2$ on ${\cal F}_1$ implies that
\begin{eqnarray}\label{eMPKK}
[H_1]|_{{\cal P}_1}=[H_2]|_{{\cal P}_1}
\end{eqnarray}
Without loss of generality, to simplify notation, we may assume
that ${\cal F}_1\supset {\cal F}$ and $\dt_2<\dt_1.$
Set ${\cal F}_2={\cal G}_1\cup {\cal F}_1.$
Let $\eta_2>0,$ $\sigma>0,$ $\dt_3>0$, $\gamma_1>0$ and a finite subset
${\cal G}_2$ be required in \ref{L2dig} corresponding to
$\ep/16,$ ${\cal F}_2,$  $L.$
We now let ${\cal G}={\cal F}_2\cup {\cal G}_2$
and $\eta=\min\{\eta_1,\eta_2\}.$  We may assume that $\dt_3<\dt_2.$
Let $x_1,x_2,...,x_m\in X$ be an
$\eta/2$-dense subset.
Suppose that
$$
 O_i\cap
O_j=\emptyset,\,\,\,{\rm if}\,\,\, i\not=j,
$$
where $O_i=\{x\in X: {\rm dist}(x, x_i)<\eta/2s\},$ and
suppose that
$$
\mu_{\tau\circ \phi}(O_i),\,\,\mu_{\tau\circ\psi}(O_i)>\eta\cdot \sigma\rforal
\tau\in T(A),1\le i\le m.
$$
We also assume that
\begin{eqnarray}\label{eMPe1}
|\tau\circ \phi(f)-\tau\circ \psi(f)|<\gamma/2\rforal f\in {\cal
G} \andeqn \rforal \tau\in T(A).
\end{eqnarray}

Since $TR(A)=0,$ there exists a
sequence of finite dimensional \SCA s $B_n$ with $e_n=1_{B_n}$ and
a sequence of \morp s $\phi_n: A\to B_n$ such that

(1) $\lim_{n\to\infty}\|e_na-ae_n\|= 0$ for all $a\in A,$

(2) $\lim_{n\to\infty}\|\phi_n(a)-e_nae_n\|=0$ for all $a\in A$
and $\phi_n(1)=e_n;$

(3) $\lim_{n\to\infty}\|\phi_n(ab)-\phi_n(a)\phi_n(b)\|=0$
for all $a, b\in A$ and

(4) $\tau(1-e_n)\to 0$ uniformly on $T(A).$

In what follows, we may assume that
\begin{eqnarray}\label{eMPT}
\tau(1-e_n)<{\sigma\cdot \eta\over{8mL+1}}\rforal \tau\in T(A).
\end{eqnarray}
 We write $B_n=\oplus_{i=1}^{r(n)}D(i,n),$ where
each $D(i,n)$ is a simple finite dimensional \CA, a full matrix
algebra. Denote by $\Phi(i,n): A\to D(i,n)$ the map which is the
composition of the projection map from $B_n$ onto $D(i,n)$ with
$\phi_n.$ Denote by $\tau(i,n)$ the standard
normalized trace on $D(i,n).$ Note that any weak limit of
$\tau(i,n)\circ \Phi(i,n)$ gives a tracial state of $A.$ By
(\ref{eMPe1}) and (4),
we have,
for all sufficiently large $n,$
$$
|\tau(i,n)\circ \Phi(i,n)\circ \phi(g)-\tau(i,n)\circ
\Phi(i,n)\circ \psi(g)| <\gamma\rforal g\in {\cal G}.
$$
Put $\phi_{(i,n)}=\Phi(i,n)\circ \phi$ and
$\psi_{(i,n)}=\Phi(i,n)\circ \psi.$ By (4), for all
sufficiently large $n,$ we have
\begin{eqnarray}\label{eMP0b1}
\mu_{\tau(i,n)\circ \phi_{(i,n)}}(O_k)\ge \sigma\cdot\eta/2
\andeqn \mu_{\tau(i,n)\circ\psi_{(i,n)}}(O_k)\ge \sigma\cdot
\eta/2,\,\,\,k=1,2,...,m
\end{eqnarray}
for each $i.$
From (1), (2), (3) above and (\ref{eMP0b1}), by applying \ref{L2dig}, for
each $i$ and all sufficiently large $n,$ we have

\begin{eqnarray}\label{eMP1}
\|\phi_{(i,n)}(f)-[(1_{D(i,n)}-p_{i,n})\phi_{(i,n)}(f)(1_{D(i,n)}-p_{i,n})+
\sum_{j=1}^mf(x_j)p(j,i,n)]\|<\ep/8,
\end{eqnarray}

\begin{eqnarray}\label{eMP2}
\|\psi_{(i,n)}(f)-[(1_{D(i,n)}-q_{i,n})\psi_{(i,n)}(f)(1_{D(i,n)}-q_{i,n})
+\sum_{j=1}^mf(x_j)q(j,i,n)]\|<\ep/8,
\end{eqnarray}
$$
\|p_{i,n}\phi_{i,n}(f)-\phi_{i,n}(f)p_{i,n}\|<\ep/8\andeqn
\|q_n\psi_{i,n}(f)-\psi_{i,n}n(f)q_{i,n}\|<\ep/8
$$
for all $f\in {\cal G},$ where $p(1,i,n),p(2,i,n),...,p(m,i,n)\in
D(i,n)$ are nonzero mutually orthogonal projections,
$p_{i,n}=1-\sum_{j=1}^mp(j,i,n)$ and
$q_{i,n}=1-\sum_{j=1}^mq(j,i,n).$ Moreover,
\begin{eqnarray}\label{eMP3}
(3mL+3)\tau(i,n)(1_{D(i,n)}-p_{i,n})<\tau(i,n)(p(k,i,n)),
\end{eqnarray}
\begin{eqnarray}\label{eMP3+}
(3mL+3)\tau(i,n)(1_{D(i,n)}-q_{i,n})<\tau(i,n)(q(k,i,n)),
\end{eqnarray}
\begin{eqnarray}\label{eMP3++}
|\tau(i,n)(p(k,i,n))-\tau(i,n)(q(k,i,n))|<
\min_k\{\tau(i,n)(p(k,i,n)\}/3mL \andeqn
\end{eqnarray}
\begin{eqnarray}\label{eMP33}
\tau(i,n)(p(k,i,n)),\,\, \tau(i,n)(q(k,i,n)>{\sigma\cdot \eta}/2,\, 1\le k\le m.
\end{eqnarray}

For  convenience, we may assume that
$$
\tau(i,n)(p(k,i,n))\ge \tau(i,n)(q(k,i,n)),\,\,\,1\le k\le m_1<m
\andeqn
$$
$$
 \tau(i,n)(q(k,i,n))\ge \tau(i,n)(p(k,i,n))\,\,\,\,\,\, m_1+1\le k\le m.
$$
In $D(i,n),$ there are projections $p(k,i,n)'\le p(k,i,n),$
$k=1,2,...,m_1$ and $q(k,i,n)'\le q(k,i,n)$ $k=m_1+1,m_1+2,...,m$
such that
$$
\tau(i,n)(p(k,i,n)')=\tau(i,n)(q(k,i,n)), k=1,2,...,m_1\andeqn
$$
$$
\tau(i,n)(q(k,i,n)')=\tau(i,n)(p(k,i,n)), k=m_1+1,m_1+2,...,m.
$$
Put $p_{i,n}'=\sum_{k=1}^{m_1}p(k,i,n)'+\sum_{k=m_1+1}^mp(k,i,n)$ and
$q_{i,n}'=\sum_{k=1}^{m_1}q(k,i,n)+\sum_{k=m_1+1}^m q(k,i,n)'.$ One
computes, by (\ref{eMP3}), (\ref{eMP3+}) and (\ref{eMP3++}), that
$$
(2mL+1)\tau(i,n)(1_{D(i,n)}-p_{i,n}')<\tau(i,n)(p(k,i,n)),\,\,\,
\tau(i,n)(p(k,i,n)') \andeqn
$$
$$
(2mL+1)\tau(i,n)(1_{D(i,n)}-q_{i,n}')<\tau(i,n)(q(k,i,n)),\,\,\,
\tau(i,n)(q(k,i,n)') \andeqn
$$
$$
\tau(i,n)(p(k,i,n)')>(\sigma\cdot\eta)/2\,\,\,{\rm for}\,\,\,
1\le k\le m_1,\andeqn
$$
$$
\tau(i,n)(q(k,i,n)>(\sigma\cdot \eta)/2\,\,\,{\rm for}\,\,\,
m_1<k\le m.
$$
One also has
$$
\|\phi_{(i,n)}(f)-[(1_{D(i,n)}-p_{i,n}')\phi_{(i,n)}(f)(1_{D(i,n)}-p_{i,n}')+
\sum_{j=1}^{m_1}f(x_j)p(j,i,n)'
+\sum_{j=m_1+1}^mf(x_j)p(j,i,n)]\|<\ep/8
$$
and
$$
\|\psi_{(i,n)}(f)-[(1_{D(i,n)}-q_{i,n}')\psi_{(i,n)}(f)(1_{D(i,n)}-q_{i,n}')+
\sum_{j=1}^{m_1}f(x_j)q(j,i,n) +\sum_{j=m_1+1}^m
f(x_j)q(j,i,n)']\|<\ep/8
$$
for all $f\in {\cal G}.$

Note that in $D(i,n),$ there is a unitary
$u_{(i,n)}$ such that
$$
u_{(i,n)}^*p(j,i,n)'u_{(i,n)}=q(j,i,n), 1\le j\le m_1\andeqn
u_{(i,n)}^*p(j,i,n)u_{(i,n)}=q(j,i,n)', m_1<j\le m.
$$
Thus without loss of generality, we may assume that
$p(j,i,n)'=q(j,i,n)$ and $p(j,i,n)'=q(j,i,n).$ Furthermore, by
changing notation if necessary, we may assume that (\ref{eMP1}) and
(\ref{eMP2}) hold as well as
\begin{eqnarray}\label{eMP4}
(2mL+1)\tau(i,n)(1-p_n)<\tau(i,n)(p(j,i,n))
\end{eqnarray}
and $p_{i,n}=q_{i,n},$ $p(j,i,n)=q(j,i,n).$ By combining all $i,$ without
loss generality, we may write that

\begin{eqnarray}\label{eMP5}
\|\phi_n(f)-[(1_{B_n}-P_n)\phi_n(f)(1_{B_n}-P_n)+\sum_{k=1}^m
f(x_k)P(k,n)]\|<\ep/4 \andeqn
\end{eqnarray}

\begin{eqnarray}\label{eMP5+}
\|\psi_n(f)-[(1_{B_n}-P_n)\psi_n(f)(1|_{B_n}-P_n)+\sum_{k=1}^mf(x_k)P(k,n)]\|<\ep/4
\end{eqnarray}
for all $f\in {\cal F}.$  Here $P(k,n)=\sum_{i=1}^{r(n)}p(k,i,n)$
and $P_n=1_{B_n}-\sum_{k=1}^mP(k,n).$

Furthermore, we also have
\begin{eqnarray}\label{eMP6}
(2mL+1)[1_{B_n}-P_n]\le  [P(k,n)]\,\,\,{\rm in}\,\,\, B_n \andeqn
t(P(k,n))\ge {\sigma\cdot \eta\over{4}} \rforal t\in T(A).
\end{eqnarray}
It follows from (\ref{eMP6}) and (\ref{eMPT}) that ($TR(A)=0$)
\begin{eqnarray}\label{eMP7}
(mL+1)[1_A-P_n]\le [P(k,n)]\,\,\, {\rm in}\,\,\, A,
\,\,\,\,k=1,2,...,m.
\end{eqnarray}
Put
$$
 H_1(f)=(1_A-P_n)h_1(f)(1_A-P_n)\andeqn
H_2(f)=(1_A-P_n)h_2(f)(1_A-P_n)
$$
for $f\in C(X).$
It follows from (1), (2) above and by (\ref{eMP5})
and (\ref{eMP5+}), for all sufficiently large $n,$  one estimates
that
\begin{eqnarray}\label{eMP10}
\phi\sim_{\ep/4} H_1\oplus (\oplus _{k=1}^mf(x_i)P(k,n)) \andeqn
\end{eqnarray}
\begin{eqnarray}\label{eMP11}
\psi\sim_{\ep/4} H_2\oplus (\oplus _{k=1}^mf(x_i)P(k,n))
\end{eqnarray}
on ${\cal G}.$ By the choice of ${\cal F}_1$  and $\dt_2,$ one
obtains that
$$
[H_1\oplus \oplus _{k=1}^mf(x_i)P(k,n)]|_{{\cal P}_1}
=[H_2\oplus \oplus _{k=1}^mf(x_i)P(k,n)]|_{{\cal P}_1}
$$
It follows that (by working in each group $K_i(A\otimes C_n)$)
$$
[H_1]|_{{\cal P}_1}=[H_2]|_{{\cal P}_1}.
$$
Denote $E=1_A-P_n$ and define $H_0': C(X)\to M_{m}(EAE)$ by
$$
H_0'(f)={\rm diag}(f(x_1),f(x_2),...,f(x_m))\rforal f\in C(X)
$$
and define $H_0=H_0'\oplus H_0'\oplus\cdots \oplus H_0': C(X)\to
M_{mN}(EAE). $ Then, by the choice of $N,$ $\eta_1,$ $\dt_1,$
${\cal G}_1$ and ${\cal P}_1$ and applying \ref{Luniq}, we obtain
a unitary $u\in M_{mN+1}(EAE)$ such that
\begin{eqnarray}\label{eMP8-}
H_1\oplus H_0 \aeps{u}_{\ep/2} H_2\oplus H_0\,\,\,{\rm on}\,\,\,
{\cal F}.
\end{eqnarray}
Rewrite $H_0(f)=\sum_{k=1}^m f(x_k)E'_k,$ where
$E_1',E_2',...,E_m'$ are mutually orthogonal projections, for
$f\in C(X).$
 By (\ref{eMP7}), there is a
unitary $W\in A$ such that $W^*E_k'W\le P(k,n),$ $k=1,2,...,m.$
Put $Q_k=P(k,n)-W^*E_k'W$ and define
$H_{00}(f)=\sum_{k=1}^mf(x_k)Q_k$ for $f\in C(X).$ Then one has
\begin{eqnarray}\label{eMP8}
H_1\oplus H_0\oplus H_{00}\sim H_1\oplus
(\oplus_{k=1}^mf(x_k)P(k,n))\andeqn
\end{eqnarray}
\begin{eqnarray}\label{eMP9}
H_2\oplus H_0\oplus H_{00} \sim H_2\oplus
(\oplus_{k=1}^mf(x_i)P(k,n))
\end{eqnarray}
on ${\cal F}.$
Finally,  by (\ref{eMP8-}), (\ref{eMP8}),
(\ref{eMP9}), (\ref{eMP10}) and (\ref{eMP11}), one obtains
$$
\phi\sim_{\ep} \psi\,\,\,\,{\rm on}\,\,\,{\cal F}.
$$

\end{proof}

Now we are ready to prove Theorem \ref{MT02}.

 {\bf Proof of Theorem
\ref{MT02}}

\begin{proof}
Let
$\eta>0$ and  $x_1,x_2,...,x_m\in X$ be an $\eta/4$-dense subset.
Choose an integer $s>0$ such that
$$
 O_i\cap
O_j=\emptyset,\,\,\,{\rm if}\,\,\, i\not=j,
$$
where $O_i=\{x\in X: {\rm dist}(x, x_i)<\eta/2s\}.$
 Let $g_i\in
C(X)$ such that $0\le g_i\le 1,$ $g_i(x)=1$ if ${\rm dist}(x,
x_i)<\eta/8s$ and $g_i(x)=0$ if ${\rm dist}(x,x_i)\ge \eta/4s,$
$i=1,2,...,m.$

Since $A$ is simple and $h_1$ is a monomorphism,
$$
\inf\{\tau(g_i):\tau\in T(A)\}=a_i>0, i=1,2,...,m.
$$
Thus $\mu_{\tau\circ h_1}(O_i)>a_i$ for all $\tau\in T(A),$
$i=1,2,...,m.$ Choose $d=\inf\{a_i: 1\le i\le m\}/\eta.$ Then
$$
\mu_{\tau\circ h_1}(O_i)>d\cdot\eta, i=1,2,...,m.
$$
Choose $\sigma=\min\{d/2, 1/2s\}.$
Let ${\cal G}={\cal
G}\cup\{g_i: 1\le i\le m\}.$ By choosing small $\gamma,$ if
$$
|\tau\circ h_1(g)-\tau\circ h_2(g)|<\gamma\,\,\,\rforal g\in {\cal
G},
$$
one also has
$$
\mu_{\tau\circ h_2}(O_i)>\sigma\circ \eta, \,\,\,i=1,2,...,m.
$$
We see then \ref{MT02} follows from \ref{MT01}.
\end{proof}

%

{\bf Proof of Theorem \ref{MT1}}.

\begin{proof}
It is an immediate consequence of Theorem \ref{MT02}.
\end{proof}

\begin{Def}\label{IDmup}
{\rm Let $X$ be a compact metric space and $A$ be a stably finite \CA.
Let $C=PM_k(C(X))P.$ Suppose that $h: C\to A$ is a unital \hm\,
and $\tau\in T(A).$ Define $\phi: C(X)\to C$ by $\phi(f)=f\cdot P$
for $f\in C(X).$ Define ${\tilde \tau}=\tau\circ \phi: C(X)\to
{\mathbb C}.$ We use ${\tilde \mu}_{\tau}$ for the probability
measure induced by ${\tilde \tau}.$ This notation will be used
below and in the proof of \ref{VL1}. Note also if $X$ is connected
and has finite dimension, then  $C$ is a full hereditary \SCA\, of
$M_k(C(X))$ consequently  $K_i(C)=K_i(C(X))$ ($i=0,1$).
}
\end{Def}

\begin{Cor}\label{MT1C1}
Let $X$ be a finite dimensional compact metric space and
let $C=PM_k(C(X))P,$
where $P\in M_k(C(X))$ is a projection.
Let $\ep>0,$  ${\cal F}\subset C$ be a finite subset.
Let $\eta>0$ be such that
$|f(x)-f(x')|<\ep/8.$
Then,
for any integer $s>0,$ any finite $\eta/2$-dense subset $\{x_1,x_2,...,x_m\}\subset X$
such that $O_i\cap O_j=\emptyset$ if $i\not=j,$ where
$$
O_j=\{x\in X: {\rm dist}(x, x_i)<\eta/2s\}
$$
and any  $1/2s>\sigma>0,$ there exist $\gamma>0,$
$\dt>0,$ a finite subset ${\cal G}\subset C$ and a finite subset
${\cal P}\subset {\bf P}(C)$ satisfying the following:

For any unital simple \CA\, $A$ with $TR(A)=0,$ any unital
${\cal G}$-$\dt$-multiplicative \morp s, $\phi, \psi: C\to A$
with
$$
|\tau\circ \phi(a)-\tau\circ\psi(a)|<\gamma\,\,\,for \,\,\,
all\,\,\, g\in {\cal G}\,\,\, and \,\,\, for\,\,\, all\,\,\,
\tau\in T(A),
$$
if
$$
{\tilde \mu}_{\tau\circ \phi}(O_i),\,\,\, {\tilde \mu}_{\tau\circ
\psi}(O_i)
>{\sigma\cdot \eta}
$$
for all $i$ and $\tau\in T(A),$ and if
$$
[\phi]|_{{\cal P}}=[\psi]|_{{\cal P}},
$$
then there exists a unitary in $A$ such that
$$
\phi\aueeps{u}_{\ep} \psi\,\,\,\,\,\,on \,\,\, {\cal F}.
$$
\end{Cor}


\begin{proof}
First let us assume that $C=M_k(C(X)).$ Let $\{e_{ij}\}$ be a
system of matrix unitis. Suppose that $\{\phi_n\}: C\to A$ be a
unital sequentially asymptotic morphism. There is a sequence of
projections $p_n\in A$  such that
$$
\lim_{n\to\infty}\|\phi_n(e_{11})-p_n\|=0.
$$
One obtains a sequence of elements $a_n\in A$ such that
$a_n\phi_n(e_{11})a_n=p_n$ and
\begin{eqnarray}\label{e4MC1}
\lim_{n\to\infty}\|a_n-p_n\|=0
\end{eqnarray}
Let $\phi_n'(c)=a_n\phi_n(c)a_n$ for $c\in e_{11}Ce_{11}$ ($\cong
C(X)$). Then $\phi_n'$ is a completely positive linear map. Since
$\phi_n'(1_C)=p_n,$ it is a \morp. Since $\sum_{i=1}^ke_{ii}=1_C,$
it is easy to check that, for all large $n,$ there are elements
$\{a_{ij}^{(n)}\}\subset A$ such that $\{a_{ij}^{(n)}\}$ forms a
system of matrix unit (with size $k$) such that
$a_{11}^{(n)}=e_n$ and $\sum_{i=1}^k e_{ii}=1_A.$ In other words,
we may write $A=M_k(e_nAe_n).$

Define $\phi_n'': C\to  M_k(e_nAe_n)$ by $\phi_n''=\phi_n'\otimes
{\rm id}_{M_k}.$  By (\ref{e4MC1}), we have
$$
\lim_{n\to\infty}\|\phi_n-\phi_n'\|=0.
$$
It follows that
$$
\lim_{n\to\infty}\|\phi_n-\phi_n''\|=0.
$$
Therefore we may assume that $\phi(e_{11})$ is a projection. Note
$\tau(\phi(f\cdot e_{11}))=(1/k){\tilde \tau}\circ \phi(f\cdot
P))$ for $f\in C(X).$ It is then clear, by identifying
$e_{11}Ce_{11}$ with $C(X),$ that we reduce the case that
$C=M_k(C(X))$ to the case that $C=C(X)$ which has been proved in
\ref{MT01}.

Now we consider the general case.  Suppose that $C=PM_k(C(X))P$
and ${\rm dim}X=d.$ It follows from 8.12 of \cite{Hu} (see also
6.10.3 of \cite{Bl}) that there exists an integer $K\ge 1$ ($K\le
2dk$) such that there is a projection $Q\in M_K(C)$ and partial
isometry $z\in M_k(C)$ such that $z^*z=P$ and $zz^*\le Q$ and
$QM_K(C)Q\cong M_l(C(X))$ for some integer $Kk\ge l>0.$ Suppose
that $\{\phi_n\}$ and $ \{\psi_n\}$ are  two unital sequentially
asymptotic  morphisms from $C$ into $A.$ Define
$\Phi_n=\phi_n\otimes {\rm id}_{M_K}$ and
$\Psi_n=\psi_n\otimes{\rm id}_{M_K}.$ Then $\{\Phi_n\}$ and
$\{\Psi_n\}$ are two unital sequentially asymptotic morphisms from
$M_K(C)$ into $M_K(A).$ Using $\Phi_n$ and $\Psi_n$ again for the
restriction of $\Phi_n$ and $\Psi_n$ on $QM_K(C)Q.$ There exists
projections $e_n, e_n'\in M_K(A)$ such that
$$
\lim_{n\to\infty}\|\Phi_n(Q)-e_n\|=0\andeqn
\lim_{n\to\infty}\|\Psi_n(Q)-e_n'\|=0
$$
One obtains $a_n, a_n'\in M_K(A)$ such that
$$
a_n\Phi_n(Q)a_n=e_n, a_n'\Psi_n(Q)a_n'=e_n\andeqn
\lim_{n\to\infty}\|a_n-e_n\|=0
$$
Thus, by replacing $\Phi_n$ by $a_n\Phi_na_n$ and $\Psi_n$ by
$a_n'\Psi a_n',$ we may assume that $\Phi_n$ and $\Psi_n$ map into
$e_nAe_n$ and $e_n'Ae_n'$ respectively. By the assumption, we may
also assume that $e_n$ and $e_n'$ are unitarily equivalent.
Without loss of generality, therefore, we may assume that
$e_n=e_n'.$ Note that we have prove the case that $C=M_k(C(X)).$
Since $C$ is a \SCA\, of $QM_K(C)Q\cong M_l(C(X)),$ one easily
concludes that this corollary holds.
\end{proof}

\section{Approximately conjugacy}

\begin{Prop}\label{IIPDM1}
Let $(X, \alpha)$ and $(X,\beta)$ be minimal
systems such that $TR(A_{\alpha})=TR(A_{\beta})=0.$
 Then $(X, \alpha)$ and $(X, \beta)$  are $C^*$-strongly
approximately flip conjugate if  and only if there exist an
isomorphism $\phi: A_{\alpha} \to A_{\beta},$ sequences of
unitaries $\{u_n\}\subset A_{\beta},$ $\{v_n\}\subset A_{\alpha}$
and  sequences of isomorphisms $\chi_n,\,\lambda_n: C(X)\to C(X)$
such that
$$
\lim_{n\to\infty}\|{\rm ad}\, u_n\circ \phi\circ j_{\alpha}(f)-
j_{\beta}\circ \chi_n(f)\|=0\andeqn \lim_{n\to\infty}\|{\rm
ad}\,v_n\circ \phi^{-1}\circ j_{\beta}(f)-j_{\alpha}\circ
\lambda_n(f)\|=0
$$
for all $ f\in C(X).$

\end{Prop}

\begin{proof}
The ``if part" is obvious. Suppose that $(X, \alpha)$ and $(X,
\beta)$ are $C^*$-strongly approximately flip conjugate. Suppose
that $\phi_n: A_{\alpha}\to A_{\beta}$ are isomorphisms such that
$[\phi_n]=[\phi_1]$ in $KL(A_{\alpha}, A_{\beta})$ for all $n$ and
$$
\lim_{n\to\infty}\|\phi_n\circ j_{\alpha}(f)-j_{\beta}\circ
\chi_n(f)\|=0 \rforal f\in C(X).
$$
Let $\{{\cal F}_n\}$ be an increasing sequence of finite subsets
of $A_{\alpha}$ for which the union $\cup_n {\cal F}_n$ is dense
in $A.$ Since $[\phi_n]=[\phi_1]$ in $KL(A_{\alpha}, A_{\beta}),$
by Theorem 2.3  in \cite{Lnctaf}, there is a unitary $u_n\in A_{\beta}$ such
that
$$
{\rm ad}\, u_n\circ \phi_1\sim_{{1/2^n}} \phi_n\,\,\,{\rm on}\,\,\,{\cal F}_n
$$
It follows immediately that
$$
\lim_{n\to\infty}\|{\rm ad}\, u_n\circ \phi_1\circ
j_{\alpha}(f)-j_{\beta}\circ \chi_n(f)\|=0 \rforal f\in C(X).
$$
The required unitaries $\{v_n\}$ can be obtained the same way.
\end{proof}

{\bf Proof of Theorem \ref{IIITM1}}

\begin{proof}
To see the ``only if " part, suppose that there exists a sequence
of isomorphisms $\phi_n : A_{\alpha}\to A_{\beta}$ and there
exists a sequence of  isomorphisms $\chi_n: C(X)\to C(X)$ such
that $[\phi_n]=[\phi_1]$ in $KL(A_{\alpha},A_{\beta})$  for all
$n$ and
\begin{eqnarray}\label{eTM1a}
\lim_{n\to\infty}\|\phi_n(j_{\alpha}(f))-j_{\beta}\circ
\chi_n(f)\|=0\rforal f\in C(X).
\end{eqnarray}
Put $\theta=[\phi_n].$ It follows from (\ref{eTM1a}) that, for any
finitely generated subgroup $G\subset \underline{K}(C(X)),$
$$
[j_{\alpha}]\times \theta|_G=[j_{\beta}\circ \chi_n]|_G
$$
for all sufficiently large $n.$ Since $[\phi_n]=\theta,$ it
follows that
$$
(\phi_n)_{\rho}=\theta_{\rho}.
$$
Therefore
$$
\lim_{n\to\infty}\| \rho_{A_{\beta}}\circ j_{\beta}\circ
\chi_n(f)-\theta_{\rho}\circ \rho_{A_{\alpha}}\circ
j_{\alpha}(f)\|=0
$$
for all $f\in C(X)_{s.a}.$ This proves the ``only if " part.

To prove the ``if" part, we apply \ref{MT02}. Since both
$A_{\alpha}$ and $A_{\beta}$ are simple amenable separable \CA s
which satisfy the UCT, if $(X, \alpha)$ and $(X, \beta)$ satisfy
the condition of the theorem, then, by  \cite{Lnann}, there is an
isomorphism $h: A_{\alpha}\to A_{\beta}$ such that $[h]=\theta.$
Moreover, for any finite subset ${\cal G}\subset C(X)$ and any
$\gamma>0,$ there exists $N>0$ such that, for all $n\ge N,$
$$
|\tau\circ h\circ j_{\alpha}(f)-\tau\circ j_{\beta}\circ
\chi_n(f)|<\gamma \rforal f\in C(X)
$$
and all $\tau\in T(A_{\beta}).$  Since, for any finitely generated
$G\subset \underline{K}(C(X)),$
$$
[h\circ j_{\alpha}]|_G=[j_{\beta}\circ \chi_n]|_G\,\,\,\,{\rm
for\,\,\, all\,\,\,sufficiently \,\,\,} n,
$$
it follows from \ref{MT02} that there are unitaries $u_n\in
A_{\beta}$ such that
$$
\lim_{n\to\infty}\|u_n^*(h\circ j_{\alpha}(f))u_n-j_{\beta}\circ
\chi_n(f)\|=0 \rforal f\in C(X).
$$
The same are argument provides unitaries $v_n\in A_{\alpha}$ such
that
$$
\lim_{n\to\infty}\|v_n^*(h^{-1}\circ
j_{\beta}(f))v_n-j_{\alpha}\circ \lambda_n(f)\|=0 \rforal f\in
C(X).
$$
\end{proof}

\begin{Def}{\rm (3.1 of \cite{LM})}\label{IIIDwc}
{\rm Let $X$ be a compact metric space and let $\alpha, \beta: X\to X$
be minimal homeomorphisms. We say that $\alpha$ and $\beta$ are
{\it weakly approximately conjugate}, if there exist two
sequences of homeomorphisms $\sigma_n, \gamma_n: X\to X$ such that
$$
\lim_{n\to\infty}(\sigma_n\circ \alpha\circ
\sigma_n^{-1})(f)=\beta(f)\andeqn \lim_{n\to\infty}(\gamma_n\circ
\beta\circ \gamma_n^{-1})(f)=\alpha(f)
$$
for all $ f\in C(X).$ }
\end{Def}

It easy to see (as in 3.2 of \cite{LM}) that there exists a
sequentially asymptotic morphism $\{\phi_n\}: A_{\alpha}\to
A_{\beta}$ and a sequentially asymptotic morphism $\{\psi_n\}:
A_{\beta}\to A_{\alpha}$ such that
$$
\lim_{n\to\infty}\|\phi_n(u_{\alpha})-u_{\beta}\|=0,\,\,\,
\lim_{n\to\infty}\|\psi_n(u_{\beta})-u_{\alpha}\|=0,
$$
\begin{eqnarray}\label{eAppd}
\lim_{n\to\infty}\|\phi_n(j_{\alpha}(f))-j_{\beta}(f\circ
\gamma_n)\|=0 \andeqn
\lim_{n\to\infty}\|\psi_n(j_{\beta}(f))-j_{\alpha}(f\circ
\sigma_n)\|=0
\end{eqnarray}
for all $f\in C(X).$

It is proved in \cite{LM} that, when $X$ is the Cantor set
$\alpha$ and $\beta$ are weakly approximate conjugate if $\alpha$
and $\beta$ have the same period spectrum. Therefore it is a
rather weaker relation.


\begin{Def}\label{IIIDconjK}
{\rm Let $X$ be a compact metric space and let $\alpha, \beta:
X\to X$ be minimal homeomorphisms. We say $\alpha$ and $\beta$ are
{\it approximately  $K$-conjugate}, if $\alpha$ and $\beta$ are
weakly approximately conjugate with the conjugate maps
$\{\sigma_n\}$ and $\{\gamma_n\}$ such that there are isomorphisms
$\phi_n: A_{\alpha}\to A_{\beta}$ and $\psi_n: A_{\beta}\to
A_{\alpha}$
$$
\lim_{n\to\infty}\|j_{\beta}(f\circ  \gamma_n)-\phi_n\circ
j_{\alpha}(f)\|=0\andeqn \lim_{n\to\infty}\|j_{\alpha}(f\circ
\sigma_n)-\psi_n\circ
j_{\beta}(f)\|=0
$$
for all $f\in C(X)$ and $[\phi_1]=[\phi_n]=[\psi_n^{-1}]$ in
$KL(A_{\alpha}, A_{\beta}).$}
This definition is only for the cases that
$TR(A_{\alpha})=TR(A_{\beta})=0.$
\end{Def}

\begin{Remark}\label{IIIMTR}
{\rm It is easy to see that approximate $K$-conjugacy is stronger
than $C^*$-strong approximately flip conjugacy. For weakly
approximate conjugacy, we do not require anything for the
conjugate maps. For approximate $K$-conjugacy, we require that the
conjugate maps have consistent information on $K$-theory. Suppose
that $\alpha$ and $\beta$ are weakly approximately conjugate with
the conjugate maps $\{\sigma_n\}$ and $\{\gamma_n\}.$ Suppose that
$\{\phi_n\}: A_{\alpha}\to A_{\beta}$ and $\{\psi_n\}:
A_{\beta}\to A_{\alpha}$ are the induced sequential asymptotic
morphisms as in \ref{IIIDwc}. Suppose that $\{\phi_n\}$ and
$\{\psi_n\}$ induce two elements $\theta\in KL(A_{\alpha},
A_{\beta})$ and in $\zeta\in KL(A_{\beta}, A_{\alpha})$  which
give unital order isomorphisms in $Hom(K_*(A_{\alpha}),
K_*(A_{\beta})$ and $Hom(K_*(A_{\beta}), K_*(A_{\alpha}),$
respectively. Furthermore, we assume that $\theta\times
\zeta=[{\rm id}_{A_{\alpha}}]$ and $\zeta\times \theta=[{\rm
id}_{A_{\beta}}].$
 Let $\chi_n(f)=f\circ \sigma_n$ and $\lambda_n(f)=f\circ \gamma_n.$
 Then it
is straight forward to check that, for any finitely generated
subgroup $G\subset \underline{K}(C(X)),$ $[j_{\alpha}]\times
\theta|_G= [j_{\beta}\circ \sigma_n^*]|_G$ and  $[j_{\beta}]\times
\theta^{-1}|_G=[j_{\alpha}\circ \lambda_n]|_G$ for all
sufficiently large $n.$ From (\ref{eAppd}), we also have
$$
\lim_{n\to\infty}\|\theta_{\rho}\circ \rho_{A_{\alpha}}\circ
j_{\alpha}(f)-\rho_{A_{\beta}}\circ j_{\beta}\circ
\lambda_n(f)\|=0\,\,\,{\rm and}\,\,
\lim_{n\to\infty}\|\theta_{\rho}^{-1}\circ \rho_{A_{\beta}}\circ
j_{\beta}(f)-\rho_{A_{\beta}}\circ j_{\alpha}\circ \chi_n(f)\|=0
$$
for all $f\in C(X)_{s.a.}.$ It follows from \ref{IIITM1} that $(X,
\alpha)$ and $(X,\beta)$ are $C^*$-strongly approximately flip
conjugate. Since they are also weakly approximate conjugate with
the conjugate maps $\{\sigma_n\}$ and $\{\gamma_n\},$ one checks
that they are in fact approximately $K$-conjugate. }
\end{Remark}

{\bf Proof of Theorem \ref{TCant}}

\begin{proof}

Most of the proof was given in \cite{LM}. In particular, the
equivalence of (iii), (iv) and (v) are given there. The equivalence
of (iii) and (vi) was given in \cite{GPS}. That (ii) implies (i)
follows from \ref{IIIMTR}. Moreover, it is obvious that (i)
implies (iv). It follows from Lemma 3 of \cite{T} that, if $\sigma\in
[[\beta]],$ then there is a unitary $u\in A_{\beta}$ such that
$u^*j_{\beta}(f)u=j_{\beta}(f\circ \sigma))$ for all $f\in C(X).$
Then an easy computation of ordered $K$-theory of $A_{\alpha}$ and
$A_{\beta}$ shows that (v) implies (ii).
\end{proof}

\section{The Rokhlin property}

\begin{Lem}\label{VL001}
Let $A$ be a unital simple separable \CA\, with stable rank one and
real rank zero and let
$\alpha\in Aut(A)$ such that $\alpha_{*0}|_G={\rm id}_G$ for some
subgroup  $G\subset K_0(A)$ for which $\rho_A(G)=\rho_A(K_0(A)).$
Then $\tau(a)=\tau\circ \alpha(a)$ for all
$\tau\in T(A)$ and for all $a\in A.$
\end{Lem}

\begin{proof}
First we claim that $\alpha_{*0}({\rm ker}\rho_A)={\rm
ker}\rho_A.$ Since $\alpha$ is an automorphism, it suffices to
show that $\alpha_{*0}({\rm ker}\rho_A)\subset {\rm ker}\rho_A.$
Let $z\in {\rm ker}\rho_A.$ Take $x\in K_0(A)_+\setminus \{0\}.$
Then $x\pm nz\ge 0$ for all positive integer $n.$
 It follows that
$\alpha_{*0}(x\pm n z)\ge 0.$ Thus $\rho_A(\alpha_{*0}(x\pm
nz))\ge 0.$ It is then easy to see that
$\rho_A(\alpha_{*0}(z))=0.$

For any projection $p\in A,$ there is $x\in G$ such that $x-[p]\in
{\rm ker}\rho_A.$ It follows that $x\ge 0.$ Since $A$ has stable rank one
 there
is a projection $q\in A$ such that $[q]=x.$ Thus $\tau(p)=\tau(q)$
for all $\tau\in T(A).$ Since $[q]\in G,$ by the above assumptions,
$[\alpha(q)]=[q].$ Again since $A$ has stable rank one,
 there is a partial isometry
$v\in A$ such that
$$
v^*v=q\andeqn vv^*=\alpha(q).
$$
By the first part of the proof, $\alpha_{*0}([p]-[q])\in {\rm
ker}\rho_A.$ In particular, $\tau(\alpha(p))=\tau(\alpha(q))$ for
all $\tau\in T(A).$ It follows that
\begin{eqnarray}\label{eV0021}
\tau(p)=\tau(q)=\tau(\alpha(q))=\tau(\alpha(p)).
\end{eqnarray}
In other words, $[p]-[\alpha(p)]\in {\rm ker}\rho_A.$

 Suppose that $a=\sum_{i=1}^n
\lambda_ip_i,$ where $\lambda_i$ are scalars and $p_1,p_2,...,p_n$
are mutually orthogonal projections. Then $\alpha(a)=\sum_{i=1}^n
\lambda_i\alpha(p_i).$ It follows from (\ref{eV0021}) that
$\tau(a)=\tau\circ \alpha(a)$ for all $\tau\in T(A).$
 Since
$A$ has real rank zero, it follows from
\cite{BP} that every self-adjoint element is a norm-limit of
self-adjoint elements with finite spectrum. It follows that
$\tau(a)=\tau\circ \alpha(a)$ for all self-adjoint elements. The
lemma then follows.
\end{proof}

\begin{Lem}\label{VL002}
Let $A$ be a unital simple separable \CA\, with $TR(A)=0$ and let
$\alpha\in Aut(A)$ such that $\alpha_{*0}|_G={\rm id}_G$ for some
subgroup $G$ of $K_0(A)$ for which $\rho_A(G)=\rho_A(K_0(A)).$
Suppose that $\{p_j\}$ is a central sequence of projections such
that $[p_j]\in G$  and define $\phi_j(a)=p_jap_j$ and
$\psi_n(a)=\alpha(p_j)a\alpha(p_j),$ $j=1,2,....$ Then
$\{\phi_n\}$ and $\{\psi_n\}$ are two sequentially asymptotic
morphisms. Suppose also that there are finite dimensional \SCA s
$B_j$ and $C_j=\alpha(B_j)$ with $1_{B_j}=p_j$ and
$1_{C_j}=\alpha(p_j)$ such that $[p_{j,i}]\in G$ for  each minimal
central projection $p_{j,i}$ of $B_j$ ($1\le i\le k(j)$) and there
are sequentially asymptotic  morphisms $\{\phi_j'\}$ and
$\{\psi_j'\}$ such that
$$
\phi_j'(a)\subset B_j,\,\,\,\,\,\, \psi_j'(a)\subset C_j,
$$
$$
\lim_{j\to\infty}\|\phi_j(a)-\phi_j'(a)\|=0 \andeqn
\lim_{j\to\infty}\|\psi_j(a)-\psi_j'(a)\|=0 \rforal a\in A.
$$
Then, for any $\ep>0$ and for any finite subset ${\cal G}\subset
A$ and a finite subset of projections ${\cal P}_0\subset M_k(A)$
for which $[p]\in G$ for all $p\in {\cal P}_0,$ there exists an
integer $J>0$ such that
$$
|\tau\circ \phi_j(a)-\tau\circ \psi_j(a)|<\ep/\tau(p_j) \rforal
a\in {\cal G}
$$
and for all $\tau\in T(A),$  and, for all $j>J,$
$$
[\phi_j(p)]=[\psi_j(p)]\,\,\,\text{in}\,\,\, K_0(A).
$$
\end{Lem}

\begin{proof}
Let $q\in M_k(A)$ be a projection such that $[q]\in G.$ Put
$\beta=\alpha\otimes {\rm id}_{M_k}.$ By replacing $A$ by
$M_k(A),$ $\alpha$ by $\beta,$ $\phi_j$ by $\phi_j\otimes {\rm
id}_{M_k}$ and $\psi_j$ by $\psi\otimes {\rm id}_{M_k},$ respectively,
 to
simplify notation, we may assume that $q\in A.$ By the assumption,
there is a partial isometry $v\in A$ such that $vv^*=q$ and
$vv^*=\alpha^{-1}(q).$ Define $B_{j,i}=p_{j,i}B_j,$
$i=1,2,...,k(j).$ Keep in mind that $B_{j,i}$ is a simple finite
dimensional \SCA.   Define
$$
\phi_{j,i}=p_{j,i}\phi_jp_{j,i},\,\,\,
\phi_{j,i}'=p_{j,i}\phi_j',\,\,\,
\psi_{j,i}=\alpha(p_{j,i})\psi_j\alpha(p_{j,i})\andeqn
\psi_{j,i}'=\alpha(p_{j,i})\psi_j',
$$
$i=1,2,...,k(j).$ By replacing $\psi_j$ by ${\rm ad}\, w_j\circ
\psi_j$ for some suitable unitaries $w_j,$ if necessary, we may
assume that $B_j$ is orthogonal to $C_j.$ There is a partial
isometry $z_{j,i}\in A$ such that
$$
z_{j,i}^*z_{j,i}=p_{j,i}\andeqn z_{j,i}z_{j,i}^*=\alpha(p_{j,i}).
$$
Then $z_{j,i}, B_{j,i}, C_{j,i}$ generate a \SCA\, $D_{j,i}\cong
M_2(B_{j,i})$ which is a simple finite dimensional \CA.

There are projections $e_{j,i}\in B_{j,i}$ and $e_{j,i}'\in
C_{j,i}$ such that
$$
\lim_{j\to\infty}\|e_{j,i}-\phi_{j,i}(q)\|=0\andeqn
\lim_{j\to\infty}\|e_{j,i}'-\psi_{j,i}(q)\|=0.
$$
Moreover,
$$
\lim_{j\to\infty}\|p_{j,i}v^*p_{j,i}vp_{j,i}-e_{j,i}\|=0\andeqn
\lim_{j\to\infty}\|p_{j,i}vp_{j,i}v^*p_{j,i}-\psi_{j,i}(\alpha^{-1}(q))\|=0.
$$
Thus
$$
[e_{j,i}]=[\psi_{j,i}(\alpha^{-1}(q))]
$$
for all large $j.$ On the other hand,
$$
\alpha(p_{j,i})q\alpha(p_{j,i})=\alpha(p_{j,i}\alpha^{-1}(q)p_{j,i}).
$$
Therefore
$$
[e_{j,i}']=[\alpha(e_{j,i})]
$$
for all large $j.$ But $[\alpha(e_{j,i})]-[e_{j,i}]\in {\rm
ker}\rho_A.$  So, for any $\tau\in T(A),$
$$
\tau(e_{j,i})=\tau(\alpha(e_{j,i})) \rforal \,\,\,{\rm
large}\,\,\,j.
$$
However, for each trace $\tau\in T(A),$ there is $\lambda_{j,i}>0$
such that $\tau|_{D_{j,i}}=\lambda_{j,i}Tr,$ where $Tr$ is the
standard trace on $D_{j,i}.$ Hence
$$
Tr(e_{j,i})=Tr(e_{j,i}') \rforal \,\,\,{\rm large}\,\,\,j.
$$
It follows that
$$
[e_{j,i}]=[e_{j,i}'] \,\,\,\,{\rm in}\,\,\, K_0(D_{j,i}) \rforal
\,\,\,{\rm large}\,\,\,j.
$$
Thus
\begin{eqnarray}\label{eVL002-1}
[\phi_j(q)]=[\psi_j(q)] \,\,\,\,{\rm in}\,\,\,K_0(A)  \rforal
\,\,\,{\rm large}\,\,\,j.
\end{eqnarray}
Since $[p_j]\in G,$ there is a unitary $Z_j$ such that
$$
Z_j^*\alpha(p_j)Z_j=p_j.
$$
Suppose that $a=\sum_{k=1}^m\lambda_ke_k,$ where $\lambda_k$ are
scalars and $e_1,e_2,...,e_m$ are mutually orthogonal projections.
Then, since $A$ has stable rank one (since $TR(A)=0$), by
(\ref{eVL002-1}), there is a  sequence of unitaries  $\{U_j\}$ in
$p_jAp_j$  such that
$$
\lim_{j\to\infty}\|U_j^*\phi_j(a)U_j-{\rm ad}\, Z_j\circ
\psi_j(a)\|=0.
$$
It follows that there is an integer $J_0>0$ such that
$$
|\tau(\phi_j(a))-\tau(\psi_j(a))|<\ep
$$
for all $\tau\in T(p_jAp_j)$ and  for all $j\ge J.$ Since
$TR(A)=0,$ the set of  self-adjoint elements with finite spectrum
is dense in $A_{s.a.}.$  The lemma follows.

\end{proof}

\begin{Lem}\label{VL1}
Let $A$ be a unital separable simple amenable \CA\, with $TR(A)=0$
satisfying the UCT, and let $\alpha\in
Aut(A)$ be such that $\alpha_{*0}|_G={\rm id}_{G}$ for some
subgroup $G$ of $K_0(A)$ for which $\rho_A(G)=\rho_A(K_0(A)).$
Suppose also that $\{p_j(l)\},$ $l=0,1,2,...,L,$ are central
sequences of projections in $A$ such that
$$
p_j(l)p_j(l')=0\,\,\,\text{if}\,\,\,l\not=l'\andeqn
\lim_{j\to\infty}\|p_j(l)-\alpha^l(p_j(0))\|=0,\,\,\, 1\le l\le L
$$
Then there exist central sequences of projections $\{q_j(l)\}$ and
 central sequences of partial isometries $\{u_j(l)\}$ such
that $q_j(l)\le p_j(l)$
$$
u_j(l)^*u_j(l)=q_j(l),\,\,\,u_j(l)u_j^*(l)=\alpha^l(q_j(0))
$$
for all large $j,$ and
$$
\lim_{j\to\infty}\|\alpha^l(q_j(0))-q_j(l)\|=0 \andeqn
\lim_{j\to\infty}\tau(p_j(l)-q_j(l))=0
$$
uniformly on $T(A).$
\end{Lem}

\begin{proof}
It follows from Theorem 5.2 of \cite{Lnduke} and Theorem 4.18 of
\cite{EG} that we may assume that $A=\overline{\cup_n A_n},$
where each $A_n$ has the form $P_nM_{k(n)}(C(X))P_n,$ where  each
$X$ is a finite dimensional compact metric space  and $P_n\in
M_{k(n)}(C(X))$ is a projection.
we may further assume that $\phi_n$ are injective.
Fix a finite
subset ${\cal F}\subset A$ and $\ep>0.$
 Let
${\cal F}_1=\cup_{l=0}^L\alpha^{-l}({\cal F}).$ Without loss of
generality, we may assume that ${\cal F}_1\subset C,$ where $C=P_n
M_{k(n)}(C(X))P_n.$

Note that $\{\alpha(p_j)\}$ is also a central sequence. Since
$\lim_{j\to\infty}\|p_j(l)-\alpha^l(p_j(0))\|=0,$
there are unitaries $w_j\in A$ such that
\begin{eqnarray}\label{eVL1-1}
\lim_{j\to\infty}\|w_j-1\|=0\andeqn w_j^*\alpha^l(p_j(0))w_j=p_j(l).
\end{eqnarray}
Let $\beta_j={\rm ad}w_j\circ \alpha.$

Let ${\cal P}\subset {\bf P}(C)$ be a finite subset. Choose an
integer $k_0>0$ such that
\begin{eqnarray}\label{eVL1-2}
{\cal P}\cap K_i(C,\Z/j\Z)=\{0\}, \,\,\,j\ge k_0.
\end{eqnarray}

Fix $\eta>0,$ and let $\{x_1,x_2,...,x_K\}$ be an $\eta/2$-net in
$X.$ Suppose that $s>0$ is an integer so that
$O_i\cap O_j=\emptyset,$ if $i\not=j,$ $i,j,=1,2,...,K,$
where
$$
O_i=\{x\in X: {\rm dist}(x, x_i)<\eta/2s\}, \,\,i=1,2,...,K.
$$
Let $f_i$ be a nonnegative continuous function in $C(X)$ such that
$0\le f_i\le 1,$ $f_i(y)=1$ if ${\rm dist}(y, x_i)\le \ep/4$ and
$f_i(y)=0$ if ${\rm dist}(y, x_i)\ge \eta/2.$ Let $a_i=f_iP_n\in
C.$ Since $A$ is simple, there are $b_{i,k}\in A$ such that
$$
\sum_{k=1}^{m(i)} b_{ik}^*a_ib_{ik}=1_A.
$$
Let $c=\max\{b_{ik}^*b_{ik}: 1\le k\le K\}$ and $M_0=\max\{m(k):
k=1,2,...,K\}.$ Put $\sigma=1/4cM_0\eta.$

Since $\{p_j(l)\}$ is central,
\begin{eqnarray}\label{eVL1-3}
\lim_{j\to\infty}\|\sum_{k=1}^{m(k)}p_j(l)b_{ik}^*p_j(l)a_ip_j(l)b_{ik}p_j(l)-
p_j(l)\|=0
\end{eqnarray}
$i=1,2,...,K.$

Since $TR(A)=0,$ there is a central  sequence
of projections $\{e_n\}$ and a sequence of finite dimensional
\SCA\, $B_n\subset A$ with $1_{B_n}=e_n,$ such that

(i) $e_nxe_n\subset_{\ep_n} B_n$ for all $x\in A,$ where $\ep_n>0$ and
$\sum_{n=1}^{\infty}\ep_n<\infty,$
and

(ii) $\lim_{n\to\infty}\tau(1-e_n)=0$ uniformly on $T(A).$

There is a subsequence $\{n(j)\}$ and projections $Q_j(l)\le
p_j(l)$ ($l=0,1,...,L$) such that
$$
\lim_{j\to\infty}\|Q_j(l)-p_j(l)e_{j(n)}\|=0.
$$
Note also that, by (i),
$$
\lim_{j\to\infty} {\rm dist}(p_je_{n(j)}, B_{n(j)})=0.
$$
There are also  sequences of projections $E_{n(j),l}\in B_{n(j)}$
such that $E_{n(j),l}E_{n(j),l'}=0,$ if $l\not=l',$
$$
\lim_{j\to\infty}\|E_{n(j),l}-Q_j(l)\|=0.
$$
Let $z_n$ be a sequence of unitaries such that
$$
\lim_{n\to\infty}\|z_j-1\|=0,\, z_j^*E_{n(j),l}z_j=Q_j(l)\andeqn
z_j^*E_{n(j),l}'z_j=Q_j(l), 0\le l\le L.
$$
Set $D_{n(j)}=Q_j(0)z_n^*B_{n(j)}z_nQ_j(0),$ $n=1,2,....$ Note
$D_{n(j)}$ is of finite dimensional. Write
$D_{n(j)}=\oplus_t^{M(j,n)}D_{n(j),t},$ where each $D_{n(j),t}$ is
simple and has rank $R(n(j),t).$  Let $d_{n(j),t}$ be a minimal
projection in $D_{n(j),t}.$ Since $TR(A)=0,$ by Theorem 7.1 of
\cite{Lntr0}, $A$ has real rank zero, stable rank one and weakly
unperforated $K_0(A).$ It follows from  \cite{BH} that
$\rho_A(K_0(A))$ is dense in $Aff(T(A)).$ Consequently, one has
$(k_0)!$ many mutually orthogonal and mutually equivalent projections
$d_{n(j), t,s}\le d_{n(j),t}$ ($s=1,2,...,(k_0)!$) such that
$[d_{n(j), t,s}]\in G$ and
\begin{eqnarray}\label{eVL1-t}
\tau(d_{n(j),t}-\sum_{s=1}^{(k_0)!}d_{n(j),t,s})<
{\tau(d_{n(j),t})\over{2^jM(j,n)}}\rforal \tau\in T(A).
\end{eqnarray}
Put
$$
g_{n(j),s,t}={\rm
diag}(\overbrace{d_{n(j),t,s},d_{n(j),t,s},...,d_{n(j),t,s}}
^{R(n(j),t)}),
$$
$$
g_{n(j),s}=\oplus_{t=1}^{M(j,n)}g_{n(j),s,t}\andeqn
g_{n(j)}=\sum_{s=1}^{(k_0)!}g_{n(j),s}.
$$
Note that $g_{n(j),s}$ commutes with every element in $D_{n(j)}.$
Put $C_{n(j),s}=g_{n(j),s}D_{n(j)}$ and
$C_{n(j)}=g_{n(j)}D_{n(j)}.$ Then $C_{n(j),s}$ and $C_{n(j)}$ are
finite dimensional \SCA s and the image of each minimal central
projection of $C_{n(j),s}$ and $C_{n(j)}$ in $K_0(A)$ are in $G.$

Combining (i) above, we note  that $\{g_{n(j),s}\}$ and
$\{g_{n(j)}\}$ are central (as $j\to \infty$) and

(1) ${\rm dist}(g_{n(j),s}xg_{n(j),s}, C_{n(j),s})\to 0$ and ${\rm
dist}(g_{n(j)}xg_{n(j)}, C_{n(j)})\to 0$ as $j\to\infty,$

(2) $\tau(Q_j-g_{n(j)})\to 0$ as $j\to\infty$ (by (\ref{eVL1-t})).

For each $x\in {\cal F}$ and $l=0,1,...,L,$
 there are $y_j(l)\in C_{n(j)}$ such that
$$
\|g_{n(j)}\alpha^{-l}(x)g_{n(j)}-y_j(l)\|\to \,\,\,{\rm
as}\,\,\,j\to\infty.
$$
Thus
$$
\|\beta_j^l(g_{n(j)})x\beta_j^l(g_{n(j)})-\beta_j^l(y_j(l))\|\to
0\,\,\,{\rm as}\,\,\, j\to\infty.
$$
Therefore
$$
{\rm dist}(\beta_j^l(d_{n(j)})x\beta_j^l(d_{n(j)}),
\beta_j^l(C_{n(j)}))\to 0,\,\,\,l=1,2,...,L
$$
as $j\to\infty.$ Define
$$
\Phi_{j,s}'(x)=g_{n(j),s}xg_{n(j),s},\,\,\,
\Phi_j'(x)=g_{n(j)}xg_{n(j)},
$$
$$
\Psi_{j,s,l}'(x)=\beta_j^l(g_{n(j),s})x\beta_j^l(g_{n(j),s})\andeqn
\Psi_{j,l}'(x)=\beta_j^l(g_{n(j)})x\beta_j^l(g_{n(j)})
$$
for $x\in A.$ Let
$$
L_{j,s,0}: g_{n(j),s}Ag_{n(j),s}\to C_{n(j),s},\,\,\,L_j:
g_{n(j)}Ag_{n(j)}\to C_{n(j)},
 $$
 $$
L_{j,s,l}: \beta_j^l(g_{n(j),s})A\beta_j^l(g_{n(j),s})\to
\beta_j^l(C_{n(j),s}) \andeqn L_{j,l}:
\beta_j^l(g_{n(j)})A\beta_j^l(g_{n(j)})\to \beta_j^l(C_{n(j)})
$$
be
 \morp s which are extensions of ${\rm id}_{C_{n(j),s}},$ ${\rm
id}_{C_{n(j)}},$ ${\rm id}_{\beta_j^l(C_{n(j),s})}$ and ${\rm
id}_{\beta_j^l(C_{n(j)})},$ $l=1,2,...,L,$ respectively. Put
$$
\Phi_{j,s}=L_{j,s,0}\circ \Phi_{j,s},\,\,\,\Phi_j=L_j\circ \Phi_j,
\,\,\,\Psi_{j,s,l}=L_{j,s,l}\circ \Psi_{j,s,l}\andeqn
\Psi_{j,l}=L_{j,l}\circ \Psi_{j,l}'.
$$ Note that $\{\Phi_{j,s}\},$ $\{\Phi_j\},$
$\{\Psi_{j,s,l}\}$ and $\{\Psi_{j,l}\}$ are sequentially
asymptotic morphisms ($1\le s\le (k_0)!$). We also have
\begin{eqnarray}\label{eVL1-k}
\Phi_j=\oplus_s^{(k_0)!}\Phi_{j,s}\andeqn
\Psi_{j,l}=\oplus_s^{(k_0)!}\Psi_{j,s,l},\,\,\,l=1,2,...,L.
\end{eqnarray}
 Let
$H_j=\imath_j\circ \Phi_j\circ \imath,$ where $\imath: C\to A$ and
$\imath_j: C_{n(j)}\to g_{n(j)}Ag_{n(j)}$ are embeddings (we may
also omit $\imath$ and $\imath_j$ when there will be no
confusion).  There is a unitary $Z_{j,l}\in A$ such that
$Z_{j,l}^*g_{n(j)}Z_{j,l}=\beta_j^l(g_{n(j)}),$ $l=1,2,...,L.$
Define $H_{j,l}={\rm ad}\, Z_{j,l}\circ \imath_{j,l}\circ
\Psi_{j,l}\circ \imath,$ where $\imath_{j,l}:
\beta_j^l(C_{j(n})\to \beta_j^l(g_{n(j)})A\beta_j^l(g_{n(j)})$ is
an embedding.

 Therefore, for all
sufficiently large $j,$ by (\ref{eVL1-k}),
\begin{eqnarray}\label{eVL1-Kk}
[H_j]|_{{\cal P}\cap K_i(C(X),\Z/j\Z)}=0\andeqn [H_j']|_{{\cal
P}\cap K_i(C(X),\Z/j\Z)}=0,\,\,\,i=0,1,0<j\le k_0.
\end{eqnarray}
Since both $H_j$ and $H_{j,l}$ factor through a finite dimensional
\SCA\, one has
\begin{eqnarray}\label{eVL1-Ktor}
[H_j]|_{{\cal P}\cap K_1(C(X)}=0\andeqn [H_{j,l}]|_{{\cal P}\cap
K_1(C(X))}=0.
\end{eqnarray}
By applying   \ref{VL002}, one computes that, for all sufficiently
large $j,$
\begin{eqnarray}\label{eVL1-K0}
[H_j]|_{{\cal P}\cap K_0(C(X))}=[H_{j,l}]|_{{\cal P}\cap
K_0(C(X))}.
\end{eqnarray}
Combining (\ref{eVL1-K0}), (\ref{eVL1-k}), (\ref{eVL1-Ktor}) and
(\ref{eVL1-2}), one has, if $j$ is sufficiently large,
\begin{eqnarray}\label{eVL1-KK}
[H_j]|_{{\cal P}}=[H_{j,l}]|_{\cal P}.
\end{eqnarray}
It follows from \ref{VL002} that
\begin{eqnarray}\label{eVL1-t1}
\lim_{j\to\infty}(\sup\{|\tau\circ\Phi_j(a)-\tau\circ {\rm ad}\,
Z_j\circ\Psi_{j,l}(a)|: \tau\in T(g_{n(j)}Ag_{n(j)})\}) =0
\end{eqnarray}
for all $a\in A.$ Hence
$$
\lim_{j\to\infty}(\sup\{|\tau\circ H_j(x)-\tau\circ
H_{j,l}(x)|:\tau\in T(g_{n(j)}Ag_{n(j)})\})=0 \rforal x\in C.
$$
It follows from (\ref{eVL1-3}) that
\begin{eqnarray}\label{eVL1-t2}
\lim_{j\to\infty}\|\sum_{k=1}^{m(i)}g_{n(j)}b_{ik}^*g_{n(j)}a_ig_{n(j)}b_{ik}^*g_{n(j)}-g_{n(j)}\|=0
\andeqn
\end{eqnarray}
\begin{eqnarray}\label{eVL1-t3}
\lim_{j\to\infty}\|\sum_{k=1}^{m(i)}\beta_j^l(g_{n(j)})b_{ik}^*\beta_j^l(g_{n(j)})
a_i\beta_j^l(g_{n(j)})b_{ik}\beta_j^l(g_{n(j)})-\beta_j^l(g_{n(j)})\|=0
\end{eqnarray}
It follows from (\ref{eVL1-t3}), when $j$ is sufficiently large,
\begin{eqnarray}\label{eVL1-t4}
\mu_{\tau}(O_k)\ge \sigma\eta/2 \rforal \tau\in
T(g_{n(j)}Ag_{n(j)}),\,\,\, k=1,2,...,K.
\end{eqnarray}
Therefore, by \ref{MT1C1}, there is unitary $v_{j,l}\in
g_{n(j)}Ag_{n(j)}$ such that
\begin{eqnarray}\label{eVL1-eq}
\lim_{j\to\infty}\|{\rm ad}v_{j,l}\circ H_{j,l}(x)-H_j(x)\|=0\rforal
x\in {\cal F}.
\end{eqnarray}
Define $z_{j,l}=g_{n(j)}v_{j,l}^*Z_{j,l}^*\beta_j^l(g_{n(j)}).$
Then
\begin{eqnarray}\label{eVL1-c1}
z_{j,l}(H_j(x)+\Psi_{j,l}(x))= z_{j,l}\Psi_{j,l}(x)\andeqn
(H_j(x)+\Psi_{j,l}(x))z_j=H_j(x)z_j.
\end{eqnarray}
Thus, by (\ref{eVL1-eq}),
 $$
\begin{array}{ll}\label{eVL1-c2}
\lim_{j\to\infty}\|z_{j,l}\Psi_{j,l}(x)-H_j(x)z_{j,l}\|&=
\lim_{j\to\infty}\|z_{j,l}\Psi_{j,l}(x)-
H_j(x)v_{j,l}^*Z_{j,l}^*\beta_j^l(g_{n(j)})\|\\\nonumber
&=\lim_{j\to\infty}\|z_{j,l}\Psi_{j,l}(x)-
v_{j,l}^*H_{j,l}(x)Z_{j,l}^*\beta_j^l(g_{n(j)})\|\\\nonumber
&=\lim_{j\to\infty}\|z_{j,l}\Psi_{j,l}(x)-
v_{j,l}^*Z_{j,l}^*(Z_{j,l}H_{j,l}(x)Z_{j,l}^*\beta_j^l(g_{n(j)}))\|\\\nonumber
&=\lim_{j\to\infty}\|z_{j,l}\Psi_{j,l}(x)-v_{j,l}^*Z_{j,l}^*
\Psi_{j,l}(x)\|=0
\end{array}
$$
for all $x\in {\cal F}.$ In other words,
\begin{eqnarray}\label{eVL1-c3}
\lim_{j\to\infty}\|z_{j,l}\Psi_{j,l}(x)-H_j(x)z_{j,l}\|=0\rforal
x\in {\cal F}.
\end{eqnarray}

Note that, for all $x\in A,$
$$
z_{j,l}x=z_{j,l}\beta_j^l(g_{n(j)})x\andeqn
xz_{j,l}=xg_{n(j)}z_{j,l}.
$$
Therefore, since $\{g_{n(j)}\}$ and $\{\beta_j^l(g_{n(j)})\}$ are
central, we have, for $x\in {\cal F},$
$$
\begin{array}{ll}
\lim_{j\to\infty}\|z_{j,l}x-xz_{j,l}\|&=\lim_{j\to\infty}\|
z_{j,l}\beta_j^l(g_{n(j)})x-
xg_{n(j)}z_{j,l}\|\\
&=\lim_{j\to\infty}\|z_{j,l}\beta_j^l(g_{n(j)})x\beta_j^l(g_{n(j)})-g_{n(j)}xg_{n(j)}z_{j,l}\|\\
&= \lim_{j\to\infty} \|z_{j,l}\Psi_{j,l}(x)-H_j(x)z_{j,l})\|=0
\end{array}
$$
On the other hand, we have
\begin{eqnarray}\label{eVL1-v1}
z_{j,l}^*z_{j,l}=g_{n(j)}\andeqn
z_{j,l}z_{j,l}^*=\beta_j^l(g_{n(j)})
\end{eqnarray}
Recall that $\beta_j^l(g_{n(j)})=w_j^*\alpha^l(g_{n(j)})w_j.$
Define $u_j(l)=w_j^*z_{j,l}.$ Then
\begin{eqnarray}\label{eVL1-v2}
(u_j(l))^*u_j(l)=z_{j,l}^*w_jw_j^*z_{j,l}=g_{n(j)}\andeqn
u_j(l)(u_j(l))^*=\alpha^l(g_{n(j)}).
\end{eqnarray}
Denote $q_j(0)=g_{n(j)}$ and $q_j(l)=\beta_j^l(g_{n(j)}),$
$l=1,2,...,L.$
 We also have, by (i) and (2) above,
\begin{eqnarray}\label{eVL1-v3}
q_j(l)\le Q_j(l)\le p_j(l)\andeqn \tau(p_j(l)-q_j(l))\to
0\,\,\,{\rm as}\,\,\, j\to\infty\,\,\,{\rm
uniformly\,\,\,on}\,\,\,T(A).
\end{eqnarray}
Finally, since (\ref{eVL1-1}),
\begin{eqnarray}\label{eVL1-v4}
\lim_{j\to\infty}\|u_jx-xu_j\|=0\rforal x\in {\cal F}.
\end{eqnarray}

\end{proof}

The following lemma is taken from \cite{LO} that has its origin
in \cite{Ks1}.

\begin{Lem}\label{LOL}
Let $A$ be a unital separable simple \CA\, with $TR(A)=0$ and
with unique tracial state
 and $\alpha\in \Aut(A)$ such that
$\alpha^r_{*0}|G={\rm id}_{G}$ for some subgroup $G\subset K_0(A)$
for which $\rho_A(G)=\rho_A(K_0(A))$ and for some integer $r\ge
1.$ Let $m\in {\mathbb N},$ $m_0\ge m$ be the smallest integer
such that $m_0 = 0 \ {\rm mod}\, r$ and $l = m+(r-1)(m_0+1).$

Suppose that $\{e_i^{(n)}\},$ $i = 0,1,...,l,$ $n = 1, 2,...,$ are
$l+1$ sequences of projections in $A$ satisfying the following:

$$
\|\alpha(e_i^{(n)})-e_{i+1}^{(n)}\|<\dt_n, \,
\lim_{n\to\infty}\dt_n=0,
$$

$$
e_i^{(n)}e_j^{(n)}=0,\, {\rm if}\,\,\, i\not=j,\,
$$

and for each $i,$ $\{e_i^{(n)}\}$ is a central sequence.
Then for each $i = 0, 1, 2, \dots, m,$
there are central sequences of projections
$\{p_i^{(n)}\}$ ($i=0,1,2,...,m$) and a  central sequence of partial
isometries $\{w_i^{(n)}\}$ such that

$$
(w_i^{(n)})^*w_
i^{(n)} = p_{i}^{(n)}\andeqn
w_i^{(n)}(w_i^{(n)})^* = p_{i+1}^{(n)},\,i=0,1,...,m - 1,
$$

where $p_i^{(n)} \le  \sum_{j=0}^{r-1}e_{i+j(m_0+1)}^{(n)}$
and
$$
\tau(\sum_{i=0}^le_i^{(n)}-\sum_{i=0}^m p_i^{(n)})\to 0
$$
as $n\to\infty$ uniformly on $T(A).$
Moreover, for each $i$,

$$
\lim_{n \rightarrow \infty}\|\alpha(p_i^{(n)}) - p_{i+1}^{(n)}\| = 0.
$$

\end{Lem}

\begin{proof}
Since $\alpha^r_{*0}|_G={\rm id}_{G},$  it follows from
\ref{VL002} that there is a  central sequence of projections
$\{q_0^{(n)}\}$ and a central sequence of partial isometries
$\{z(0,j,n)\}$ such that $q_0^{(n)}\le e_0^{(n)},$
$$
z(0,j,n)^*=q_0^{(n)},\,\, z(0,j,n)z(0,j,n)^*=\alpha^{rj}(q_0^{(n)})
$$
and
$$
\tau(e_0^{(n)}-q_0^{(n)})\to 0.
$$
Moreover, there are central sequences of projections
$\{q_{rj}^{(n)}\}$ such that $q_{rj}^{(n)}\le e_{rj}^{(n)}$ such
that
$$
\lim_{n\to\infty}\|q_{rj}^{(n)}-\alpha^{rj}(q_0^{(n)})\|=0,\,\,\,j=1,2,...,[l/r].
$$
Since
$q_0^{(n)}\le e_0^{(n)}$ and $\|\alpha(e_0^{(n)})-e_1^{(n)}\|<\dt_n,$
there exists a projection $q_1^{(n)}\le e_1^{(n)}$ such
that
$$
\|\alpha(q_0^{(n)})-q_1^{(n)}\|<2\dt_n.
$$
Similarly, we obtain projections $q_i^{(n)}\le e_i^{(n)}$ such that
$$
\lim_{n\to\infty}\|\alpha(q_{i-1}^{(n)})-q_i^{(n)}\|=0,\,\,\,
i=2,3,...,r-1.
$$
Moreover, it is easy to check that
$$
\lim_{n\to\infty}\|\alpha^i(q_0^{(n)})-q_i^{(n)}\|=0, i=1,2,...,r-1.
$$
Since $\{q_0^{(n)}\}$ is central, $\{\alpha^i(q_0^{(n)})\}$ is also central
for each $i.$
It follows that $\{q_i^{(n)}\}$ are central for $i=1,2,...,r-1$ and
for $i=rj,$ $j=1,2,....[l/r]$
Again, there are projections $q_{i+rj}^{(n)}\le e_{i+rj}^{(n)}$
such that
$$
\lim_{n\to\infty}\|\alpha^{i+rj}(q_0^{(n)})-q_{i+rj}^{(n)}\|=0,\,i=1,2,...,r-1.
$$
Define $z_{i,j,n}'=\alpha^i(z_{0,j,n}').$
Then
$$
(z_{i,j,n}')^*z_{i,j,n}'=\alpha^i(z_{0,j,n}^*z_{0,j,n})=\alpha^i(q_0^{(n)})
\andeqn
z_{i,j,n}'(z_{i,j,n}')^*=\alpha^i(z_{0,j,n}z_{0,j,n}^*)=
\alpha^{i+rj}(q_0^{(n)}).
$$
In particular, $\{z_{i,j,n}\}$ is central, since $\{z_{0,j,n}\}$
is. Because
$$
\lim_{n\to\infty}\|\alpha^i(q_0^{(n)})-q_{i}^{(n)}\|=0 \andeqn
\lim_{n\to\infty}\|\alpha^{i+rj}(q_0^{(n)})-\alpha^{rj}(q_{i}^{(n)})\|=0\,\,\,
i=1,2,...,r-1,
$$
one obtains  sequences of unitaries $\{Z_j^{(n)}\}$ and
$\{U_j^{(n)}\}$ such that
$$
\lim_{n\to\infty}\|Z_j^{(n)}-1\|=0,\,\,\,
(Z_j^{n})^*\alpha^i(q_0^{(n)})Z_j^{(n)}=q_{j}^{(n)}, \andeqn
$$
$$
\lim_{n\to\infty}\|U_j^{(n)}-1\|=0,\,\,\,
U_j^{(n)}\alpha^{i+rj}(q_0^{(n)})(U_j^{(n)})^*=\alpha^{rj}(q_{i}^{(n)})
$$
$j=1,2,...,r-1.$
Define $z_{i,j,n}=U_j^{(n)}z_{i+1,j,n}'Z_j^{(n)}.$ Then $\{z_{i, j,n}\}$ is
central and
$$
z_{i,j,n}^*z_{i,j,n}=q_{i}^{(n)}\andeqn
z_{i,j,n}z_{i,j,n}^*=\alpha^{rj}(q_{i}^{(n)}).
$$

Now put
$$
p_i^{(n)}=\sum_{j=0}^{r-1}q_{i+j(m_0+1)}^{(n)}, i=0,1,2,...,m.
$$
Note
that
$$
\tau(\sum_{i=0}^le_i^{(n)}-\sum_{i=0}^mp_i^{(n)})
<\sum_{i=0}^l\tau(e_i^{(n)}-q_i^{(n)})=l\cdot\tau(e_0^{(n)}-q_0^{(n)})\to 0
$$

To prove so defined $p_i^{(n)}$ meets the requirements, one checks
exactly the same way as in the proof of Lemma 3.2 \cite{LO}.
\end{proof}

Let $\{E_{i,j}\}$ be a system of matrix units and ${\cal K}$ be
the compact operators on $\ell^2({\mathbb Z})$ where we identify
$E_{i,i}$ with the one-dimensional projection onto the functions
supported by $\{i\} \subset {\mathbb Z}$.
Let $S$ be the canonical shift operator on $\ell^2({\mathbb Z})$.
Define an automorphism $\sigma$ of
${\cal K}$ by $\sigma(x) = SxS^*$ for all $x \in {\cal K}$.
Then $\sigma(E_{i,j}) = E_{i + 1, j + 1}$.
For any $N \in {\mathbb N}$
let $P_N = \sum_{i = 0}^{N - 1}E_{i,i}$.

To prove Theorem \ref{MTrok}, we quote the following lemma.

\begin{Lem}{\rm (Kishimoto, 2.1 of \cite{Ks1})}\label{LK}
{\rm For any $\eta > 0$ and $n \in {\mathbb N}$ there exist $N \in
{\mathbb N}$ and projections $e_0, e_1, \dots, e_{n - 1}$ in
${\cal K}$ such that

$$
\begin{array}{l}
\sum_{i=0}^{n - 1}e_i \leq P_N\\
\|\sigma(e_i) - e_{i + 1}\| < \eta, \ i = 0, \dots, n - 1, \ e_n = e_0\\
\frac{n\dim e_0}{N} > 1 - \eta.
\end{array}
$$
}
\end{Lem}

{\bf Proof of Theorem \ref{VTRC}}

\begin{proof}
We use a modified argument of Theorem 3.4 in \cite{LO} by applying
\ref{LOL}. We proceed as follows.

Let $\ep>0.$  Let $\ep/2 > \eta > 0$ and $m \in {\mathbb N}$ be
given. Choose $N$ which satisfies the conclusion of Lemma \ref{LK}
(with this $\eta$ and $n=m$). Identify $P_N{\cal K}P_N$ with
$M_N.$ Let ${\cal G}=\{E_{i+1,i}:i=0,1,...,N-1\}$ be a set of
generators of $M_N.$ Let $e_0, e_1,...,e_{m-1}$ be as in the
conclusion of Lemma \ref{LK}.
 For any $\ep > 0,$
there is $\dt > 0$ depends only on $N$
such that, if

$$
\|ag-ga\| <\dt
$$
for $g \in {\cal G}$,
then

$$
\|ae_i - e_ia\| < \ep/2, \,\,\, i = 0, 1,..., n.
$$
We assume that $\dt < \eta.$  Fix a finite subset ${\cal
F}_0\subset A.$
Choose $m_0 \in {\mathbb N}$ such that $m_0 \geq m$ is the
smallest integer with $m_0 = 0 \ {\rm mod}\ r$.
Let $L = N + (r - 1)(m_0 + 1)$.

Since $\alpha$ has the tracial Rokhlin property, there exists
a sequence of projections
$\{e_i^{(k)}: i = 0, 1,...,L\}$
satisfying the following:

$$
\|\alpha(e_i^{(k)})-e_{i+1}^{(k)}\|<{\dt\over{(2^k)4N}},\,\,\,
e_i^{(k)}e_j^{(k)}=0,\,\,\,
{\rm if}\,\,\, i \not =  j,
$$

$$
\lim_{k\to\infty}\|e_i^{(k)}a-ae_i^{(k)}\| = 0
\rforal a \in A,\,i = 0, 1,..., L \andeqn
$$

$$
\tau(1-\sum_{i = 0}^{L-1}e_i^{(k)}) < \eta/2 \rforal \tau \in
T(A),\ k = 1, 2,....
$$

By applying Lemma \ref{LOL}, we obtain a central sequence
$\{w_i^{(k)}\}$ in $A$ such that

$$
\begin{array}{l}
(w_i^{(k)})^*w_i^{(k)} = P_0^{(k)}   \andeqn\\
w_i^{(k)}(w_i^{(k)})^* = P_i^{(k)}, \ k = 0, 1,...,\
i = 0, 1, \dots, N, \\
P_i^{(k)}P_j^{(k)} = 0 \ i \not= j,\\
\|\alpha(P_i^{(k)}) - P_{i + 1}^{(k)}\| < \frac{\dt}{4L},
\ k = 0, 1,...,\ i = 0, 1, \dots, N - 1,\\
\tau(\sum_{i=0}^{L-1}e_i^{(k)}- \sum_{i=  0}^{N-1}P_i^{(k)}) <
\eta/2, \ \rforal \tau \in T(A)
\end{array}
$$
where $P_i^{(k)} \le \sum_{j = 0}^{r-1}e_{i+(m_0+1)j}^{(k)}$ for
$i = 0, 1, \dots, N$.

It follows that $\{ \alpha^l(w^{(k)} )\},$ $l = 0, 1,..., N$ are
all central sequences. As the same argument in Lemma \ref{LOL}
there is a unitary $u_k \in U(A)$ with $\|u_k - 1\|<\dt/2N$ such
that ${\rm ad}\, u_k \circ \alpha(P_i^{(k)}) = P_{i+1}^{(k)},$ $i
= 0, 1,..., N-1.$ Put $\beta_k = {\rm ad}\, u_k\circ \alpha$, and
$w^{(k)} = w_0^{(k)}$. Choose a large $k,$ such that

$$
\|\beta_k^{l}(w^{(k)})a - a\beta_k^{l}(w^{(k)})\| < \dt\rforal
a\in {\cal F}_0,
$$
$l = 0, 1, ..., N.$

Now let $C_1$ and  $C_2$
be the \CA s generated by
$w^{(k)}, \beta_k^{1}(w^{(k)}), ..., \beta_k^{N - 1}(w^{(k)})$ and by
$w^{(k)}, \beta_k^{1}(w^{(k)}), ..., \beta_k^{N}(w^{(k)}),$
respectively. Note that $C_1\cong M_{N},$ $C_2\cong M_{N + 1}.$
Define a \hm\, $\Phi: C_1\to {\cal K}$ by

$$
\Phi(\beta_k^{i}(w^{(k)})) = E_{i+1,i}, \ i = 0, 1, ..., N - 1
$$

(see Lemma \ref{LK}).
Then one has $\sigma\circ \Phi|_{C_1} = \Phi\circ
\beta_k|_{C_1}$ and $\Phi(C_1) = P_{N}{\cal K}P_{N}.$
Now we apply Lemma \ref{LK} to obtain
mutually orthogonal projections $e_0, e_1,...,e_{m-1}$
in $M_N$ such that
$$
\|\sigma(e_i)-e_{i-1}\|<\eta\andeqn
{m\dim e_0\over{N}} > 1 - \eta.
$$
Let
$p_i = \Phi^{-1}(e_i),$ $i = 0, 1,..., m - 1.$
One estimates that
$$
\tau(\sum_{i = 0}^{N-1}P_i^{(k)} - \sum_{i=0}^{m-1}p_i) <
1 - \sum_{i=0}^{m-1}\frac{\dim(e_0)}{N}
= 1 - \frac{m\dim(e_0)}{N} < \eta < \frac{\ep}{2}
$$
for all $\tau \in T(A).$
So one has mutually orthogonal projections
$p_0, p_1, p_2,..., p_{m - 1}$ such that
$$
\|\beta_k(p_i) - p_{i+1}\| < \frac{\ep}{2},
 i = 0, 1, 2,..., m - 1,\, p_m = p_0.
$$
By the choice of $\dt,$ one also has
$$
\|ap_i - p_ia\| < \ep, \  i = 0, 1,..., m - 1, \ \rforal a \in
{\cal F}_0
 \andeqn
$$

$$
\tau(1-\sum_{i=0}^{m - 1}p_i) <
\tau(1-\sum_{i=0}^{L-1}e_i^{(k)})+\tau(\sum_{i=0}^{L-1}e_i^{(k)} -
\sum_{i=0}^{N-1}P_i^{(k)}) + \frac{\ep}{2} < \eta/2+\eta/2 +
\frac{\ep}{2}< \ep,
$$
for all $\tau \in T(A).$ Since
$$
\|\beta_k - \alpha\| < \dt/2 < \ep/2,
$$
one finally has
$$
\|\alpha(p_i) - p_{i+1}\| < \ep,\,\,\,i = 0, 1,..., m - 1\,\ , p_m
= p_0.
$$
In other words, $\alpha$ has the tracial cyclic Rokhlin property.
\end{proof}

{\bf Proof of Theorem \ref{MTAH}}

\begin{proof}
This follows from Theorem \ref{MTrok} and Theorem 3.4 of
\cite{LO}.
\end{proof}

{\bf Proof of Theorem \ref{VTRC}}

\begin{proof}
The fact that (i), (ii) and (iii) are equivalent (without assuming
that $\alpha^r_{*0}|_G={\rm id}_{G}$) is established in \cite{OP}.

That (iv) $\Rightarrow$ (v) is given by Theorem 2.9 in \cite{LO}.
It is known that (v) $\Rightarrow$ (iv).

Thus we have shown that (i), (ii), (iii), (iv) and (v) are
equivalent.

 To see these imply (vi),
we apply the classification theorem in \cite{Lnduke}. It follows
that $A\rtimes_{\alpha}\Z$ is a unital simple AH-algebra with no
dimension growth and with real rank zero. By (ii), it has a unique
tracial state.

It is obvious that (vi) implies (iii).
\end{proof}

{\bf Proof of Theorem \ref{TLFur}}

\begin{proof}
It follows from Theorem 8.3 in \cite{OP}  that $\alpha_{\theta,
\gamma, d,f}$ has the tracial Rokhlin property. It also follows
from Corollary 8.4 in \cite{OP} that
$A_{\theta}\rtimes_{\alpha}\Z$  is a unital simple \CA\, with real
rank zero and with unique tracial state. It follows from Theorem
\ref{VTRC} that $A_{\theta}\rtimes_{\alpha}\Z$ is an AH-algebra
with no dimension growth and with real rank zero. This proves (1).
To see (2), since $d=0,$ by Lemma 1.7 of \cite{OP},
$K_0(A_{\theta}\rtimes_{\alpha}\Z)$ is torsion free. It follows
that it is an A$\T$-algebra.
\end{proof}

\end{document}